\documentclass{elsart1p}

\usepackage{amsfonts}
\usepackage{graphicx}
\usepackage{amsmath}
\usepackage{amssymb}

 \advance\hoffset by -0.0 truecm

\newtheorem{theorem}{Theorem}
\newtheorem{corollary}{Corollary}
\newtheorem{lemma}{Lemma}
\newtheorem{proposition}{Proposition}
\theoremstyle{definition}
\newtheorem{definition}{Definition}
\theoremstyle{remark}
\newtheorem{remark}{Remark}
\numberwithin{equation}{section}

\DeclareMathOperator{\spn}{span} \DeclareMathOperator{\Cl}{Cl}
\DeclareMathOperator{\dist}{dist} 
\DeclareMathOperator{\Spin}{Spin} \DeclareMathOperator{\SU}{SU}
\DeclareMathOperator{\SO}{SO} 

\DeclareMathOperator{\arctanh}{arctanh}
\DeclareMathOperator{\arcsinh}{arcsinh}
\DeclareMathOperator{\arccosh}{arccosh}
\DeclareMathOperator{\arccotanh}{arccotanh}
\DeclareMathOperator{\arccot}{arccot}
\DeclareMathOperator{\cotanh}{cotanh}
\DeclareMathOperator{\cotan}{cotan}

\journal{Journal des Math\'ematiques Pures et Appliqu\'ees}

\begin{document}

\begin{frontmatter}
\title{Sub-Lorentzian Geometry on Anti-de Sitter Space}

\author[Der-Chen]{Der-Chen Chang\thanksref{grant1}}
\ead{chang@georgetown.edu}
\author[AlexIrina]{Irina Markina\thanksref{grant2}}
\ead{irina.markina@uib.no}
\author[AlexIrina]{Alexander Vasil'ev\thanksref{grant2}\corauthref{cor}}
\corauth[cor]{Corresponding author.}
\ead{alexander.vasiliev@math.uib.no}

\thanks[grant1]{The first author has been supported by a research grant
from the United States Army Research Office, by a competitive
research grant of the Georgetown University, and by the grant of the Norwegian Research Council \#180275/D15.}

\thanks[grant2]{The second and the third authors have been  supported by the grants of the Norwegian Research Council \# 177355/V30, \#180275/D15, and by the European Science Foundation Networking Programme HCAA}

\address[Der-Chen]{Department of Mathematics, Georgetown University, Washington
D.C. 20057, USA}
\address[AlexIrina]{Department of Mathematics,
University of Bergen, Johannes Brunsgate 12, Bergen 5008, Norway}

\begin{abstract}
Sub-Riemannian Geometry is proved to play an important role in many applications, e.g., 
Mathematical Physics and Control Theory. 
Sub-Riemannian Geometry enjoys major differences from the Riemannian being a generalization of the latter at the same time, e.g., geodesics are not unique and may be singular, the Hausdorff dimension
is larger than the manifold topological dimension. There exists a large amount of literature developing sub-Riemannian Geometry. However, very few is known about its extension to  pseudo-Riemannian analogues. It is natural to begin such a study with some low-dimensional manifolds.
Based on ideas from sub-Riemannian geometry we develop  sub-Lorentzian geometry over the
classical 3-D anti-de Sitter space. Two different distributions of the tangent bundle of anti-de Sitter space yield two different geometries: sub-Lorentzian and sub-Riemannian. We use Lagrangian and Hamiltonian formalisms for both sub-Lorentzian and sub-Riemannian geometries to find geodesics. 

\medskip

\noindent
{\bf R\'esum\'e} 
\medskip

\noindent
Il a \'et\'e prouv\'e que la G\'eom\'etrie sub-Riemannienne joue un r\^ole important dans des nombreuses applications, par ex., Physique Math\'ematique et Th\'eorie de Contr\^ole.  La G\'eom\'etrie sub-Riemannienne a des diff\'erences consid\'erables par rapport a celle Riemannienne, \'etant au m\^eme temps une g\'en\'eralisation de celle-ci, par ex., les g\'eod\'esiques ne sont pas uniques et peuvent \^etre singuli\`eres, la dimension de Hausdorff est plus grande que la dimension topologique de vari\'et\'e. Il y a une quantit\'e importante de litt\'erature qui d\'eveloppe la G\'eom\'etrie sub-Riemannienne.  Cependant, on conna{\^\i}t tr\`es peu sur son extension naturelle aux analogues pseudo Riemanniens. C'est naturel de commencer une telle \'etude avec des vari\'et\'es de basse dimension. En se basant sur les id\'ees de la g\'eom\'etrie sub-Riemannienne, on d\'eveloppe la g\'eom\'etrie sub-Lorentzienne sur l'espace anti-de Sitter classique. Deux distributions diff\'erentes du faisceau tangent de l'espace d'anti-de Sitter donnent deux g\'eom\'etries diff\'erentes: sub-Lorentzienne et sub-Riemannienne.  On utilise \'egalement les formalismes de Lagrange et d'Hamilton pour les deux g\'eom\'etries, sub-Lorentzienne et sub-Riemanniene, pour trouver les g\'eod\'esiques.
\end{abstract}

\begin{keyword}{Sub-Riemannian and sub-Lorentzian geometries, geodesic, anti-de Sitter space, Hamiltonian system, Lagrangian, spin group, spinors}
\MSC Primary: 53C50, 53C27; Secondary: 83C65
\end{keyword}
\end{frontmatter}

\section{Introduction}

Many interesting studies of anticommutative algebras and 
sub-Riemannian structures may be seen in a general setup of
Clifford algebras and spin groups. Among others we distinguish the
following example. The unit 3-dimensional sphere $S^3$ being embedded into the Euclidean space $\mathbb R^4$ 
possesses a clear manifold structure with the Riemannian metric. It is interesting to consider the sphere $S^3$ as an algebraic object $S^3=\SO(4)/\SO(3)$ where the group $\SO(4)$ preserves the global Euclidean metric of the ambient space $\mathbb R^4$ and $\SO(3)$ preserves the Riemannian metric on $S^3$. The quotient $\SO(4)/\SO(3)$ can be realized as the group $\SU(2)$ acting on $S^3$ as on the space of  complex vectors $z_1,z_2$ of unit norm $|z_1|^2+|z_2|^2=1$. It is isomorphic to the group of  unit quaternions with the group operation given by the quaternion multiplication. It is natural to make the correspondence between $S^3$ as a smooth manifold and $S^3$ as a Lie group acting on this manifold.  
The corresponding Lie algebra is given by left-invariant vector fields with
non-vanishing commutators. This leads to  construction of
a sub-Riemannian structure on $S^3$, see \cite{CChM} (more about sub-Riemannian geometry see, for instance,~\cite{Mon,Str1,Str2,Sub}). The commutation relations for vector fields on the tangent bundle of $S^3$ come from the non-commutative multiplication for quaternions. Unit quaternions, acting by conjugation on vectors from $\mathbb R^3$ (and $\mathbb R^4$), define  rotation in $\mathbb R^3$ (and $\mathbb R^4$), thus preserving the positive-definite metric in $\mathbb R^4$. At the same
time, the Clifford algebra over the vector space $\mathbb R^3$
with the standard Euclidean metric gives rise to the spin group
$\Spin(3)=\SU(2)$ that acts on the group of unit spinors in the
same fashion leaving some positive-definite quadratic form invariant. Two models are equivalent but the latter admits
various generalizations. We are primary aimed at switching the
Euclidean world to the Lorentzian one and  sub-Riemannian geometry
to sub-Lorentzian following a simple example similar to the above of a
low-dimensional space that leads us to  sub-Lorentzian geometry
over the pseudohyperbolic space $H^{1,2}$ in $\mathbb R^{2,2}$. 
In General Relativity the simply connected covering manifold of $H^{1,2}$ 
is called the universal anti-de Sitter space~\cite{Ol,Pen,Sz}.

We start with some more rigorous explanations. A real Clifford
algebra is associated with a vector space $V$ equipped with a quadratic
form $Q(\cdot,\cdot)$. The multiplication (let us denote it by $\otimes$) in the Clifford algebra
satisfies the relation $$v\otimes v=-Q(v,v)1,$$ for $v\in V$, where 1 is the
unit element of the algebra. We restrict ourselves to $V=\mathbb
R^3$ with two different quadratic forms:
$$Q_{\mathcal E}(v,v)=\mathcal E v\cdot v,\quad \mathcal
E=\left[\array{rrr} 1 & 0 & 0
\\
0 & 1 & 0
\\
0 & 0 & 1
\endarray\right],
$$ and
$$Q(v,v)= Iv\cdot v,\quad
I=\left[\array{rrr} -1 & 0 & 0
\\
0 & 1 & 0
\\
0 & 0 & 1
\endarray\right].
$$
The first case represents the standard inner product in the Euclidean space $\mathbb
R^3$.  The second case corresponds to the Lorentzian metric in
$\mathbb R^3$ given by the diagonal metric tensor  with the
signature $(-,+,+)$. The corresponding Clifford algebras
we denote by $\Cl(0,3)=\Cl(3)$ and $\Cl(1,2)$. The basis of the
Clifford algebra $\Cl(3)$ consists of the elements
$$\{1,i_1,i_2,i_3,i_1\otimes i_2,i_1\otimes i_3,i_2\otimes
i_3,i_1\otimes i_2\otimes i_3\},\ \ \text{with} \ \ i_1\otimes
i_1=i_2\otimes i_2=i_3\otimes i_3=-1.$$ The algebra $\Cl(3)$ can
be associated with the space $\mathbb H\times \mathbb H$, where
$\mathbb H$ is the quaternion algebra. The basis of the
Clifford algebra $\Cl(1,2)$ is formed by
$$\{1,e,i_1,i_2,e\otimes i_1,e\otimes i_2,i_1\otimes
i_2,e\otimes i_1\otimes i_2\},\ \ \text{with} \ \ e\otimes e=1,\ \
i_1\otimes i_1=i_2\otimes i_2=-1.$$ In this case the algebra is
represented by $2\times 2$ complex matrices.

Spin groups are generated by quadratic elements of Clifford algebras.
We obtain the spin group $\Spin(3)$ in the case of the Clifford algebra $\Cl(3)$, and the group
$\Spin(1,2)$ in the case of the Clifford algebra
$\Cl(1,2)$. The group $\Spin(3)$ is represented by the group
$\SU(2)$ of unitary $2\times 2$ complex matrices with determinant
$1$. The elements of $\SU(2)$ can be written as $$\left[\array{rr}
a & b
\\
-\bar b & \bar a
\endarray\right], \quad a,b\in\mathbb C,\quad |a|^2+|b|^2=1.$$
The group $\Spin(3)=\SU(2)$ forms a double cover of the group of
rotations $\SO(3)$. In this case the Euclidean metric in $\mathbb
R^3$ is preserved under the actions of the group $\SO(3)$. The
group $\Spin(3)=\SU(2)$ acts on spinors similarly to how
$\SO(3)$ acts on vectors from $\mathbb R^3$. Indeed, given an
element $R\in \SO(3)$ the rotation is performed by the matrix
multiplication $RvR^{-1}$, where $v\in \mathbb R^3$.  An element
$U\in\SU(2)$ acts over spinors regarded as $2$ component vectors
$z=(z_1,z_2)$ with complex entries in the same way $UzU^{-1}$.
This operation defines a `half-rotation' and preserves the
positive-definite metric for spinors. Restricting ourselves to 
spinors of length 1, we get the manifold $\{(z_1,z_2)\in \mathbb
C^2:\ |z_1|^2+|z_2|^2=1\}$ which is the unit sphere $S^3$.

Now we turn to the Lorentzian metric and to the Clifford algebra
$\Cl(1,2)$. The spin group $\Spin^+(1,2)$ is represented by the
group $\SU^+(1,1)$ whose elements are
$$\left[\array{rr} a & b
\\
\bar b & \bar a
\endarray\right], \quad a,b\in\mathbb C,\quad |a|^2-|b|^2=1.$$
The group $\Spin^+(1,2)=\SU^+(1,1)$ forms a double cover of the
group of Lorentzian rotations $\SO(1,2)$ preserving the
Lorentzian metric $Q(v,v)$. Acting on spinors, the group
$\Spin^+(1,2)=\SU^+(1,1)$ preserves the pseudo-Riemannian metric for
spinors. Unit spinors $(z_1,z_2)$, $|z_1|^2-|z_2|^2=1$, are invariant under
the actions of the corresponding group $\Spin^+(1,2)=\SU^+(1,1)$. The
manifold $H^{1,2}=\{(z_1,z_2)\in \mathbb C^2:\ |z_1|^2-|z_2|^2=1\}$ is a
$3$-dimensional Lorentzian manifold known as a pseudohyperbolic space in Geometry and
as the anti-de Sitter
space $AdS_3$ in General Relativity. In fact, anti-de Sitter space is the maximally symmetric,
simply connected, Lorentzian manifold of constant negative
curvature. It is one of three maximally symmetric cosmological constant solutions
to Einstein's field equation: de Sitter space with a positive
cosmological constant $\Lambda$, anti-de Sitter space with a negative
cosmological constant $-\Lambda$, and the flat space. Both de
Sitter $dS_3$ and anti-de Sitter $AdS_3$ spaces may be treated as non-compact
hypersurfaces in the corresponding pseudo-Euclidean spaces $\mathbb R^{1,3}$ and $\mathbb R^{2,2}$. Sometimes de
Sitter space $dS_3$ or the hypersphere is used as a direct analogue to the sphere $S^3$ given its
positive curvature. However, $AdS_3$ geometrically is a natural object for us to work with. We reveal the
analogy between $S^3$
and $AdS_3$ as follows. The group of rotations $\SO(4)$ in the usual
Euclidean $4$-dimensional space acts as  translations on the Euclidean
sphere $S^3$ leaving it invariant. As it has been mentioned at the beginning, the sphere $S^3$ can be thought
of as the Lie group $S^3=\SO(4)/\SO(3)$ endowed with the group law
given by the multiplication of  matrices from $\SU(2)$ which is
the multiplication law for unit quaternions. The Lie algebra is
identified with the left-invariant vector fields from the tangent
space at the unity. The tangent bundle admits the natural
sub-Riemannian structure and $S^3$ can be considered as a
sub-Riemannian manifold. This geometric object was studied in
details in~\cite{CChM}. It appears throughout celestial mechanics in works of Feynman and Vernon who described it in the language of two-level systems, in Berry's phase in quantum mechanics or in the Kustaaheimo-Stifel transformation for regularizing binary collision. 

Instead of $\mathbb R^4$, we consider now the space
$$\mathbb R^{2,2}=\{\ (x_1,x_2,x_3,x_4)\in\mathbb R^4\ \ \text{with
a pseudo-metric\ } dx^2=-dx_1^2-dx_2^2+dx_3^2+dx_4^2\}.$$ The
group $\SO(2,2)$ acting on $\mathbb R^{2,2}$ is a direct analog
of the rotation group $\SO(4)$ acting on $\mathbb R^4$. We consider $AdS_3$ as a manifold $H^{1,2}=\SO(2,2)/\SO(1,2)$ with the
Lorentzian metric induced from $\mathbb R^{2,2}$. Sometimes in physics literature, $AdS_3$ appears as a universal cover of $H^{1,2}$. It is worth to mention that $H^{1,2}$ is a homogeneous non-compact manifold and the group
$\SO(2,2)$ acts as an isometry on $H^{1,2}$. The difference between this construction and above mentioned sphere is that $S^3$ itself is a group, whereas $H^{1,2}$ is not. However, $\SO(2,2)$ can be factorized  as $\SO(2,2)=\SU^+(1,1)\times \SU^+(1,1)'$, so  $H^{1,2}$ becomes a group manifold for $\SU^+(1,1)$, and topologically they are the same. The group law is defined by the matrix multiplication of  elements from $\SU^+(1,1)$. The reader can find more information about the group actions and relation to General Relativity, e.~g.~\cite{Nab,Por}. 
Left-invariant vector fields on the tangent bundle are not
commutative and this gives us an opportunity to consider an
analogue of sub-Riemannian geometry, that is called {\it
sub-Lorentzian} geometry on $\SU^+(1,1)$ (which by abuse of notation, we  call the  $AdS$ group). 
The geometry of anti-de Sitter space was studied in numerous works, see, for example,~\cite{BS,Car,Lu,Nat,ROC}.

Very few is known about extension of sub-Riemannian geometry to its pseudo-\-Rieman\-nian analogues.
The simplest example of a sub-Riemannian structure is provided by the 3-D Heisenberg group. Let us mention that recently Grochowski studied its sub-Lorenzian analogue \cite{Groch1, Groch2}.  Our approach deals with non-nilpotent groups over $S^3$ and $AdS_3$.

The paper is organized in the following way. In Section 2 we give the precise form of left-invariant vector fields defining sub-Lorentzian and sub-Riemannian structures on anti-de Sitter group. In Sections 3 and 4 the question of existence of smooth horizontal curves in the sub-Lorentzian manifold is studied. The Lagrangian and Hamiltonian formalisms are applied to find sub-Lorentzian geodesics in Sections 5 and 6. Section 7 is devoted to the study of a sub-Riemannian geometry defined by the distribution generated by spacelike vector fields of anti-de Sitter space. In both sub-Lorentzian and sub-Riemannian cases we find geodesics explicitly.

\section{Left-invariant vector fields}

We consider the  $AdS$ group topologically as a $3$-dimensional manifold  $H^{1,2}$ in
$R^{2,2}$
$$H^{1,2}=\{(x_1,x_2,x_3,x_4)\in\mathbb R^{2,2}:\ -x_1^2-x_2^2+x_3^2+x_4^2=-1\},$$ and as a group  $\SU^+(1,1)$ with the group
law  given by the multiplication of the matrices from $\SU^+(1,1)$.
We write $a=x_1+ix_2$, $b=x_3+ix_4$, where $i$ is the complex
unity. For each matrix $\left[\array{rr} a & b
\\
\bar b & \bar a
\endarray\right]\in\SU^+(1,1)$ we associate its coordinates to the
complex vector $p=(a,b)$. Then the multiplication law between $p=(a,b)$ and $q=(c,d)$ written in
coordinates is
\begin{equation}\label{ml}
pq=(a,b)(c,d)=(ac+b\bar d,ad+b\bar c).\end{equation} 
Then,  $AdS$ 
with the multiplication law~\eqref{ml} is the Lie group with the unity $(1,0)$, with the inverse to $p=(a,b)$ element $p^{-1}=(\bar a,-b)$, and with the left translation $L_p(q)=pq$. The Lie algebra is associated with the left-invariant vector fields at the identity of the group.
To calculate
the real left-invariant vector fields, we write the multiplication law~\eqref{ml}
in real coordinates, setting $c=y_1+iy_2$, $d=y_3+iy_4$. Then
\begin{equation}\label{ml1}
\begin{split} pq & = (x_1,x_2,x_3,x_4)(y_1,y_2,y_3,y_4) \\
 & = (x_1y_1-x_2y_2+x_3y_3+x_4y_4,\ x_2y_1+x_1y_2+x_4y_3-x_3y_4,\\ &\ \ \quad  
x_3y_1+x_4y_2+x_1y_3-x_2y_4,\
x_4y_1-x_3y_2+x_2y_3+x_1y_4).\end{split}\end{equation} The
tangent map $(L_p)_*$ corresponding to the left translation $L_p(q)$ is
\begin{equation*}(L_p)_*  = \left[\array{rrrr} x_1 & -x_2 & x_3 &
x_4
\\
x_2 & x_1 & x_4 & -x_3
\\
x_3 & x_4 & x_1 & -x_2
\\
x_4 & -x_3 & x_2 & x_1
\endarray\right].\end{equation*} The left-invariant vector fields
are the left translations of vectors at the unity by the tangent
map $(L_p)_*$: $\widetilde X=(L_p)_*X(0)$. Letting $X(0)$ be the
vectors of the standard basis in $\mathbb R^{2,2}$ (that
coincides with the Euclidean basis in $\mathbb R^4$), we get the
left-invariant vector fields \[\begin{array}{lll} \widetilde X_1 &
=x_1\partial_{x_1}+x_2\partial_{x_2}+x_3\partial_{x_3}+x_4\partial_{x_4},
\\ \widetilde X_2 & = -x_2\partial_{x_1}+x_1\partial_{x_2}+x_4\partial_{x_3}-x_3\partial_{x_4},
\\ \widetilde X_3 & = x_3\partial_{x_1}+x_4\partial_{x_2}+x_1\partial_{x_3}+x_2\partial_{x_4},
\\ \widetilde X_4 & = x_4\partial_{x_1}-x_3\partial_{x_2}-x_2\partial_{x_3}+x_1\partial_{x_4}
\end{array}\] in the basis
$\partial_{x_1},\partial_{x_2},\partial_{x_3},\partial_{x_4}$. Let
us introduce the matrices
\begin{equation*}U  = \left[\array{rrrr} 1 & 0 & 0 & 0
\\
0 & 1 & 0 & 0
\\
0 & 0 & 1 & 0
\\
0 & 0 & 0 & 1
\endarray\right],\qquad J   = \left[\array{rrrr} 0 & 1 & 0 & 0
\\
-1 & 0 & 0 & 0
\\
0 & 0 & 0 & -1
\\
0 & 0 & 1 & 0
\endarray\right],\end{equation*} \begin{equation*}
E_1  = \left[\array{rrrr} 0 & 0 & 1 & 0
\\
0 & 0 & 0 & 1
\\
1 & 0 & 0 & 0
\\
0 & 1 & 0 & 0
\endarray\right],\qquad E_2 = \left[\array{rrrr} 0 & 0 & 0
& 1
\\
0 & 0 & -1 & 0
\\
0 & -1 & 0 & 0
\\
1 & 0 & 0 & 0
\endarray\right].\end{equation*} Then the left-invariant vector
fields can be written in the form $$\widetilde
X_1=xU\cdot\nabla_x,\quad \widetilde X_2=xJ \cdot\nabla_x,\quad \widetilde
X_3=xE_1 \cdot\nabla_x,\quad \widetilde X_4=xE_2\cdot\nabla_x,$$ where
$x=(x_1,x_2,x_3,x_4)$,
$\nabla_x=(\partial_{x_1},\partial_{x_2},\partial_{x_3},\partial_{x_4})$ and ''$\cdot$'' is the dot-product in $\mathbb R^4$.
The matrices possess the following properties:
\begin{itemize}
\item[$\bullet$]{Anti-commutative rule or the Clifford algebra
condition:
\begin{equation}\label{m1}
J E_1 +E_1 J =0,\quad E_2E_1 +E_1 E_2=0,\quad J E_2+E_2J =0.
\end{equation}}
\item[$\bullet$]{Non-commutative rule:
\begin{equation}\label{m2}[\frac{1}{2}J ,\frac{1}{2}E_1 ]=\frac{1}{4}(J E_1 -E_1 J )=\frac{1}{2}E_2,
\quad [\frac{1}{2}E_2,\frac{1}{2}E_1 ]=\frac{1}{2}J ,\quad
[\frac{1}{2}J ,\frac{1}{2}E_2]=-\frac{1}{2}E_1 .
\end{equation}}
\item[$\bullet$]{Transpose matrices:
\begin{equation}\label{m3}J ^{T}=-J ,\quad E_2^{T}=E_2,\quad E_1 ^{T}=E_1 .\end{equation}}
\item[$\bullet$]{Square of matrices:
\begin{equation}\label{m4}J ^{2}=-U,\quad E_2^{2}=U,\quad E_1 ^{2}=U.\end{equation}}
\end{itemize} 
As a consequence we obtain
\begin{itemize}
\item[$\bullet$]{Product of matrices:
\begin{equation}\label{m5}J E_1 =E_2,\quad E_2E_1 =J ,\quad J E_2=-E_1 .\end{equation}}
\end{itemize}

The inner $\langle\cdot,\cdot\rangle$ product in $\mathbb
R^{2,2}$ is given by \begin{equation}\label{inpr}\langle  x,y\rangle=\mathcal I x\cdot y,\
\ \text{with}\ \ \mathcal I  = \left[\array{rrrr} -1 & 0 & 0 & 0
\\
0 & -1 & 0 & 0
\\
0 & 0 & 1 & 0
\\
0 & 0 & 0 & 1
\endarray\right].\end{equation} Given the inner product~\eqref{inpr} we have
\begin{equation}\label{m6}\langle x,xE_1 \rangle=\langle x,xJ \rangle=\langle x,xE_2\rangle=0,\end{equation}
\begin{equation}\label{m7}\langle xJ ,xE_1 \rangle=\langle xE_2,xE_1 \rangle=\langle xJ ,xE_2\rangle=0,\end{equation}
\begin{equation}\label{m8}\langle xJ ,xJ \rangle=-1,\quad\langle xE_2,xE_2\rangle=\langle xE_1 ,xE_1 \rangle=1.\end{equation}

The vector field $\widetilde X_1$ is orthogonal to $AdS$. Indeed, if
we write $AdS$ as a hypersurface $F(x_1,x_2,x_3,x_4)=-x_1^2-x_2^2+x_3^2+x_4^2+1=0$,
then
$$\frac{dF(c(s))}{ds}=2\Big(-x_1\frac{dx_1}{ds}-x_2\frac{dx_2}{ds}+x_3\frac{dx_3}{ds}+x_4\frac{dx_4}{ds}\Big)
=\langle
\widetilde X_1,\frac{dc(s)}{ds}\rangle=0$$ for any smooth curve
$c(s)=(x_1(s),x_2(s),x_3(s),x_4(s))$ on $AdS$. From now on we denote
the vector field $\widetilde X_1$ by $N$. Observe, that $|N|^2=\langle N,N\rangle=-1$. Up to certain ambiguity we use the same notation $|\cdot|$ as the norm (whose square is not necessary positive) of a vector and as the absolute value (non-negative) of a real/complex number.
Other vector fields are orthogonal to $N$ with respect to the
inner product $\langle\cdot,\cdot\rangle$ in $\mathbb R^{2,2}$:
$$\langle N,\widetilde X_2\rangle=\langle N,\widetilde X_3\rangle=\langle N,\widetilde
X_4\rangle=0.$$ We conclude that the vector fields $\widetilde
X_2$, $\widetilde X_3$, $\widetilde X_4$ are tangent to $AdS$.
Moreover, they are mutually orthogonal with 
$$|\widetilde X_2|^2=\langle \widetilde X_2,\widetilde X_2\rangle=-1,
\quad |\widetilde X_3|^2=|\widetilde X_4|^2=1.$$ We denote the
vector field $\widetilde X_2$ by $T$ providing time orientation
 (for the terminology see the end of the present section). The spacelike vector fields $\widetilde X_3$ and
$\widetilde X_4$ will be denoted by $X$ and $Y$ respectively. We
conclude that $T,X,Y$ is the basis of the tangent bundle of $AdS$.
In Table~\ref{t1} the commutative relations between $T,X,$ and $Y$
are presented.
\begin{table}[ht] \caption{Commutators of left-invariant vector fields} \label{t1}
\begin{center}
\begin{tabular}{|| c | r | r | r ||}
\hline \     &  $T$ & $X$ & $Y$ \\  \hline $T$ &  $  0  $ & $2Y$ &
$-2X$ \\ \hline $X$ & $-2Y$ & $ 0 $ & $-2T$
\\ \hline $Y$ & $2X$ & $2T$ & $0$
\\ \hline
\end{tabular}
\end{center}
\end{table} We see that if we fix two of the vector fields, then
they generate, together with their commutators, the tangent bundle of
the manifold $AdS$. 
\medskip 

\begin{definition}
Let $M$ be a smooth $n$-dimensional manifold, $\mathcal D$ be a smooth $k$-dimensional,
$k<n$, bracket generating distribution on $TM$, and  $\langle \cdot, \cdot \rangle_{\mathcal D}$ be a smooth Lorentzian metric on $\mathcal D$. Then the triple $(M,\mathcal D, \langle \cdot, \cdot \rangle_{\mathcal D})$ is called the {\it sub-Lorentzian manifold}.
\end{definition}

We deal with two following cases in Sections~\ref{tx}--\ref{hamilton} and Section~\ref{xy} respectively:
\begin{itemize}
\item[1.]{The horizontal distribution $\mathcal D$ is generated by
the vector fields $T$ and $X$: $\mathcal D=\spn\{T,X\}$. In this
case $T$ provides the time orientation and $X$ gives the spatial direction on $\mathcal D$.  The
 direction
$Y=\frac{1}{2}[T,X]$, orthogonal to the distribution $\mathcal D$,
is the second spatial direction on the tangent bundle. The metric $\langle \cdot, \cdot \rangle_{\mathcal D}$ is given by the restriction of $\langle \cdot, \cdot \rangle$ from $\mathbb R^{2,2}$. This case corresponds to the sub-Lorentzian manifold $(AdS,\mathcal D, \langle \cdot, \cdot \rangle_{\mathcal D})$.}
\item[2.]{The horizontal distribution $\mathcal D$ is generated by
the vector fields $X$ and $Y$: $\mathcal D=\spn\{X,Y\}$. In this
case both of the directions are spatial on $\mathcal D$. The  direction
$T=\frac{1}{2}[Y,X]$, orthogonal to the distribution $\mathcal D$. In this case, the triple  $(AdS,\mathcal D, \langle \cdot, \cdot \rangle_{\mathcal D})$ is a sub-Riemannian manifold.}
\end{itemize}

The ambient metric with the signature $(-,-,+,+)$ of $\mathbb R^{2,2}$ restricted to the tangent bundle $TAdS$ of $AdS$ is the Lorentzian
metric with the signature $(-,+,+)$, and therefore, $AdS$ is a Lorentzian manifold. The vector fields $T,X,Y$ form an orthonormal basis of each tangent space $T_pAdS$ at 
$p\in AdS$. We introduce a time orientation on $AdS$. A vector $v\in
T_pAdS$ is said to be {\it timelike} if $\langle v,v\rangle<0$,
{\it spacelike} if $\langle v,v\rangle>0$ or $v=0$, and {\it lightlike} if
$\langle v,v\rangle=0$ and $v\neq 0$. By previous consideration we have $T$ as
a timelike vector field and $X,Y$ as spacelike vector fields at each $p\in
AdS$. A timelike  vector $v\in T_pAdS$ is said to be future-directed
if $\langle v,T\rangle<0$ or past-directed if $\langle
v,T\rangle>0$. A smooth curve $\gamma: [0,1]\to AdS$ with
$\gamma(0)=p$ and $\gamma(1)=q$ is called timelike (spacelike,
lightlike) if the tangent vector $\dot\gamma(t)$ is timelike (spacelike,
lightlike) for any $t\in [0,1]$. If $\Omega_{p,q}$ is the non-empty set of
all timelike, future-directed smooth curves $\gamma(t)$ connecting
the points $p$ and $q$ on $AdS$, then the distance between $p$ and $q$ is
defined as
$$
\dist:=\sup\limits_{\gamma\in\Omega_{p,q}}\int\limits_{0}^1\sqrt{-\langle
\dot\gamma(t),\dot\gamma(t)\rangle} dt.
$$
A geodesic in any manifold $M$ is a curve $\gamma:\ [0,1]\to M$ whose vector 
field is parallel, or equivalently, geodesics are the curves of acceleration zero. 
A manifold $M$ is called geodesically connected if, given two
points $p,q\in M$, there is a geodesic curve $\gamma(t)$ connecting
them. Anti-de Sitter space $AdS$ is not geodesically connected, see
\cite{Hawking,Neill}.

The concept of causality is important in the study of Lorentz manifolds.
We say that $p\in M$ chronologically (causally) precedes $q\in M$
if there is a timelike (non-spacelike) future-directed (if
non-zero) curve starting at  $p$ and ending at $q$. For each $p\in M$ we
define the chronological future of $p$ as
$$
I^+(p)=\{q\in M:\,\, \mbox{$p$ chronologically precedes $q$}\},
$$
and the causal future of $p$ as
$$
J^+(p)=\{q\in M:\,\, \mbox{$p$ causally precedes $q$}\}.
$$
The  conformal infinity due to Penrose is timelike. One can make
analogous definitions replacing `future' by `past'.

From the mathematical point of view the spacelike curves have the same right to be studied as timelike or lightlike curves. Nevertheless, the timelike curves and lightlike curves possess an additional physical meaning as the following example shows. 

\bigskip

\noindent {\bf Example 1.} Interpreting the $x_1$-coordinate of $AdS$ as time measured in
some inertial frame ($x_1=t$), the timelike curves represent motions of particles such that 
$$\Big(\frac{dx_2}{dt}\Big)^2+\Big(\frac{dx_3}{dt}\Big)^2<1.$$ It is assumed that units have been 
chosen so that $1$ is the maximal allowed velocity for a matter particle (the speed of light). 
Therefore, timelike curves represents motions of matter particles. Timelike geodesics represent   motions with constant speed. In addition, the length $$\tau (\gamma)=\int\limits_{0}^1\sqrt{-\langle
\dot\gamma(t),\dot\gamma(t)\rangle}\,dt, $$ of a timelike curve $\gamma:\ [0,1]\to AdS$ 
is interpreted as the proper time measured by a particle between events $\gamma(0)$ and $\gamma(1)$.
\medskip

Lightlike curves represent motions at the speed of light and the lightlike geodesics represent motions along the light rays. 

\section{Horizontal curves with respect to the distribution $\mathcal D=\spn\{T,X\}$}\label{tx}

Up to Section~\ref{xy} we shall work with the horizontal distribution $\mathcal D=\spn\{T,X\}$ and the Lorentzian metric on $\mathcal D$, which is the restriction of the metric $\langle\cdot,\cdot\rangle$ from $\mathbb R^{2,2}$.
We say that an absolutely continuous curve $c(s):[0,1]\to AdS$ is
{\it horizontal} if the tangent vector $\dot c(s)$ satisfies the relation $\dot
c(s)=\alpha(s)T(c(s))+\beta(s)X(c(s))$.
\medskip 

\begin{lemma}\label{l1} A curve
$c(s)=(x_1(s),x_2(s),x_3(s),x_4(s))$ is horizontal with respect to
the distribution $\mathcal D=\spn\{T,X\}$, if and only if,
\begin{equation}\label{hc1}-x_4\dot x_1+x_3\dot x_2 -x_2\dot x_3 +x_1\dot x_4 =0\qquad \text{or}\qquad \langle x E_2,\dot c\rangle=0.\end{equation}
\end{lemma}

\begin{pf}
The tangent vector to the curve $c(s)=(x_1(s),x_2(s),x_3(s),x_4(s))$
written in the left-invariant basis  $(T,X,Y)$ admits the form $$\dot
c(s)=\alpha T+\beta X + \gamma Y.$$ Then
$$\gamma=\langle\dot c, Y\rangle=\mathcal I\dot c\cdot Y=- x_4\dot x_1+ x_3\dot x_2- x_2\dot x_3+x_1\dot x_4
=\langle xE_2,\dot c\rangle.$$ We conclude that
$$\gamma=0,$$ if and only if, the condition~\eqref{hc1} holds.\qed
\end{pf}

In other words, a curve $c(s)$ is horizontal, if and only if, its velocity
vector $\dot c(s)$ is orthogonal to the missing direction $Y$. The
left-invariant coordinates $\alpha (s)$ and $\beta(s)$ of a
horizontal curve $c(s)=(x_1(s),x_2(s),x_3(s),x_4(s))$ are
\begin{equation}\label{alpha}
\alpha=\langle\dot c, T\rangle=x_2\dot x_1 - x_1\dot x_2 +x_4\dot
x_3 - x_3\dot x_4=\langle xJ ,\dot c\rangle,
\end{equation}
\begin{equation}\label{beta}
\beta=\langle\dot c, X\rangle=-x_3\dot x_1 - x_4\dot x_2 +x_1\dot
x_3 + x_2\dot x_4=\langle xE_1 ,\dot c\rangle.
\end{equation}

Let us write the definition of the horizontal distribution
$\mathcal D=\spn\{T,X\}$ using the contact form. We define the
form $\omega=-x_4 dx_1 +x_3 dx_2- x_2 dx_3+ x_1 dx_4=\langle
xE_2,dx\rangle$. Then,
$$\omega(N)=0,\quad \omega(T)=0,\quad \omega(X)=0,\quad
\omega(Y)=1,$$  and $\ker\omega=\spn \{N,T,Y\}$, The horizontal
distribution can be defined as follows $$\mathcal D=\{V\in TAdS:\
\omega(V)=0\},\quad\text{or}\quad \mathcal D=\ker \omega\cap TAdS,
$$ where $TAdS$ is the tangent bundle of $AdS$.

The length $l(c)$ of a horizontal curve $c(s):\ [0,1]\to AdS$ is
defined by the following formula
$$l(c)=\int_{0}^{1}|\langle\dot c(s),\dot c(s)\rangle|^{1/2}\,ds.$$
Using the orthonormality of the vector fields $T$ and $X$, we
deduce that 
$$l(c)=\int_{0}^{1}\big|-\alpha^2(s)+\beta^2(s)\big|^{1/2}\,ds.$$
We see that the restriction onto
the horizontal distribution $\mathcal D\subset TAdS$ of the non-degenerate metric $\langle\cdot
,\cdot\rangle$ defined on $TAdS$ gives the Lorentzian
metric which is non-degenerate. The definitions of timelike (spacelike, lightlike) horizontal vectors $v\in\mathcal D_p$ are the same as for the vectors $v\in T_pAdS$. A horizontal curve $c(s)$
is timelike (spacelike, lightlike) if its velocity vector $\dot c(s)$ is horizontal timelike (spacelike, lightlike) vector at each point of this curve.
\medskip

\begin{lemma}
Let $\gamma(s)=(y_1(s),y_2(s),y_3(s),y_4(s))$ be a horizontal timelike future-directed (or past-directed) curve and $c(s)=L_p(\gamma(s))$ be its left translation by $p=(p_1,p_2,p_3,p_4)$, $p\in AdS$. Then the curve $c(s)$ is horizontal timelike and future-directed (or past-directed).
\end{lemma}

\begin{pf}
Let us denote by $(c_1(s),c_2(s),c_3(s),c_4(s))$ the coordinates of the curve $c(s)$. Then, by~\eqref{ml1} we have  
\begin{equation}\label{cc}
\begin{array}{lll}
c_1(s) & =p_1y_1(s)-p_2y_2(s)+p_3y_3(s)+p_4y_4(s), \\
c_2(s) & =p_2y_1(s)+p_1y_2(s)+p_4y_3(s)-p_3y_4(s), \\
c_3(s) & =p_3y_1(s)+p_4y_2(s)+p_1y_3(s)-p_2y_4(s), \\
c_4(s) & =p_4y_1(s)-p_3y_2(s)+p_2y_3(s)+p_1y_4(s).
\end{array}
\end{equation}
Differentiating with respect to $s$, we calculate the horizontality condition~\eqref{hc1} for the curve $c(s)$. Since $-p_1^2-p_2^2+p_3^2+p_4^2=-1$,  straightforward simplifications lead to the relation
\begin{equation*}
\langle\dot c, Y\rangle=-c_4\dot c_1+c_3\dot c_2-c_2\dot c_3+c_1\dot c_4  =(-p_1^2-p_2^2+p_3^2+p_4^2)(-y_4\dot y_1+y_3\dot y_2-y_2\dot y_3+y_1\dot y_4 )=0,
\end{equation*} and the curve $\gamma$ is horizontal.

Let us show that the curve $c(s)$ is timelike and future-directed provided $\gamma(s)$ is such. We calculate $$\langle\dot c, T\rangle=c_2\dot c_1 - c_1\dot c_2 +c_4\dot
c_3 - c_3\dot c_4=(p_1^2+p_2^2-p_3^2-p_4^2)(y_2\dot y_1 - y_1\dot y_2 +y_4\dot
y_3 - y_3\dot y_4)=\langle\dot\gamma, T\rangle$$ and 
$$\langle\dot c, X\rangle=-c_3\dot c_1 - c_4\dot c_2 +c_1\dot
c_3 + c_2\dot c_4=(p_1^2+p_2^2-p_3^2-p_4^2)(-y_3\dot y_1 - y_4\dot y_2 +y_1\dot
y_3 + y_2\dot y_4)=\langle\dot\gamma, X\rangle$$
from~\eqref{alpha},~\eqref{beta}, and~\eqref{cc}. Since the horizontal coordinates are not changed, we conclude that the property timelikeness and future-directness is preserved under the left translations.\qed
\end{pf}

In view that the left-invariant coordinates of the velocity vector to a horizontal curve do not change under left translations, we conclude the following analogue of the preceding lemma.
\medskip 

\begin{lemma}
Let $\gamma(s)=(y_1(s),y_2(s),y_3(s),y_4(s))$ be a horizontal spacelike (or lightlike) curve and $c(s)=L_p(\gamma(s))$ be its left translation by $p=(p_1,p_2,p_3,p_4)$, $p\in AdS$. Then the curve $c(s)$ is horizontal spacelike (or lightlike).
\end{lemma}

\section{Existence of smooth horizontal curves on $AdS$}\label{connectivity}

The question of the connectivity by geodesics of two arbitrary points on a
Lorentzian manifold is not trivial, because we have to distinguish
timelike  and spacelike curves. The problem becomes more
difficult if we study connectivity for  sub-Lorentzian geometry. In the classical Riemannian geometry all geodesics can be found as solutions to the Euler-Lagrange equations and they coincide with the solutions to the corresponding Hamiltonian system obtained by the Legendre transform. In the sub-Riemannian geometry, any solution to the Hamiltonian system is a horizontal curve and satisfies the Euler-Lagrange equations. However, a solution to the Euler-Lagrange equations is a solution to the Hamiltonian system only if it is horizontal. 

In the case of sub-Lorentzian geometry we have no information about such a correspondence. As it will be shown in Sections~6 and~7 the solutions to the Hamiltonian system are horizontal. It is a rather expectable fact given the corresponding analysis of sub-Riemannian structures, e.~g., on nilpotent groups, see~\cite{BGG,CChGr}. Since $\{T,X,Y=1/2[T,X]\}$ span the tangent space at each point of $AdS$ the existence of horizontal curves is guaranteed by Chow's theorem~\cite{Chow}. So as the first step, in this section we study connectivity by smooth horizontal curves. The main results states that any two points can be connected by a smooth horizontal curve. A naturally arisen  question is whether the found horizontal curve is timelike (spacelike, lightlike)?

First, we introduce a parametrization of $AdS$ and present the
horizontality condition and the horizontal coordinates in terms
of this parametrisation.

The manifold $AdS$ can be parametrized by
\begin{equation}\label{e2}
\begin{array}{lll}
x_1 & = & \cos a\cosh\theta,   \\ 
x_2 & = & \sin a\cosh\theta,  \\ 
x_3 & = & \cos b\sinh\theta,\\
x_4 & = & \sin b\sinh\theta,
\end{array}
\end{equation}
 with $a,b\in (-\pi,+\pi]$, $\theta\in
(-\infty,\infty)$. Setting $\psi=a-b$, $\varphi=a+b$, we formulate the
following lemma.
\medskip

\begin{lemma}\label{l2}
Let $c(s)=(\varphi(s),\psi(s),\theta(s))$ be a curve on $AdS$. The
curve is horizontal, if and only if,
\begin{equation}\label{e1}
\dot\varphi\cos\psi\sinh 2\theta-2\dot\theta\sin\psi=0.
\end{equation} The horizontal coordinates $\alpha$ and $\beta$ of the velocity vector
are \begin{equation}\label{palpha}
\alpha=-\frac{1}{2}(\dot\varphi\cosh 2\theta+\dot\psi)=-\dot a\cosh^2\theta-\dot b\sinh^2\theta,
\end{equation}
\begin{equation}\label{pbeta}
\beta=\frac{1}{2}(\dot\varphi\sin\psi\sinh 2\theta+2\dot\theta\cos\psi).
\end{equation}
\end{lemma}

\begin{pf}
Using the parametrisation~\eqref{e2} of $AdS$, we calculate
\begin{equation}\label{e3}
\begin{array}{lll}
 \dot x_1 & = & -\dot a\sin
a\cosh\theta+\dot\theta\cos a\sinh\theta, \\
\dot x_2 & = & \dot a\cos
a\cosh\theta+\dot\theta\sin
a\sinh\theta,  \\ \dot x_3 & =  &-\dot b\sin
b\sinh\theta+\dot\theta\cos b\cosh\theta,\\
\dot x_4 & =  &\dot b\cos
b\sinh\theta+\dot\theta\sin
b\cosh\theta.
\end{array} 
\end{equation}
Substituting the expressions for $x_k$ and $\dot x_k$, $k=1,2,3,4$, in~\eqref{hc1},~\eqref{alpha}, and~\eqref{beta},
in terms of $\varphi$, $\psi$ and $\theta$, we get the
necessary result.
\qed
\end{pf}

We also need the following obvious technical lemma formulated without proof.
\medskip

\begin{lemma}\label{l3}
Given $q_0$, $q_1$, $I\in\mathbb R$, there is a smooth function
$q:[0,1]\to\mathbb R$, such that $$q(0)=q_0,\qquad q(1)=q_1,\qquad
\int_{0}^{1}q(u)\,du=I.$$
\end{lemma}
\medskip

\begin{theorem}\label{t2}
Let $P$ and $Q$ be two arbitrary points in $AdS$.  Then there is a smooth horizontal curve
joining $P$ and $Q$.
\end{theorem}

\begin{pf}
Let $P=P(\varphi_0,\psi_0,\theta_0)$ and
$Q=Q(\varphi_1,\psi_1,\theta_1)$ be coordinates of the points
$P$ and~$Q$. In order to find a horizontal curve $c(s)$ we must solve equation~\eqref{e1} with the
boundary conditions \begin{eqnarray*} c(0) &
=P,\quad\text{or}\quad\varphi(0)=\varphi_0,\quad
\psi(0)=\psi_0,\quad \theta(0)=\theta_0, \\
c(1) & =Q,\quad\text{or}\quad\varphi(1)=\varphi_1,\quad
\psi(1)=\psi_1,\quad \theta(1)=\theta_1.
\end{eqnarray*}
Assume that $\sin\psi\neq 0$ we rewrite the equation~\eqref{e1}
as \begin{equation}\label{e4}
2\dot\theta=\dot\varphi\cot\psi\sinh 2\theta.
\end{equation} To simplify matters, let us introduce two new smooth functions $p(s)$ and $q(s)$ by $$2\theta(s)=\arcsinh
p(s),\quad\psi(s)=\arccot
q(s),$$ and let  the function $\varphi(s)$ is set  as $\varphi(s)=\varphi_0+s(\varphi_1-\varphi_0)$. Then we will define  the smooth functions
$p(s)$ and $q(s)$ satisfying the horizontality condition~\eqref{e4} for $c=c(s)$.
Let $k=\varphi_1-\varphi_0$. Then equation~\eqref{e4} admits the
form $$\frac{\dot p(s)}{\sqrt{1+p^2(s)}}=kp(s)q(s).$$ Separation
of variables leads to the equation $$\frac{d
p}{p\sqrt{1+p^2}}=kq(s)\,ds,$$ that after integrating
gives
$$-\arctanh\frac{1}{\sqrt{1+p^2(s)}}=k\Big(\int_{0}^{s}q(\tau)\,d\tau+C\Big).$$
To define the constant $C$, we use the boundary conditions at
$s=0$. Observe that $$\frac{1}{\sqrt{1+p^2(0)}}=\frac{1}{\cosh 2\theta_0}\quad\text{
and}\quad \frac{1}{\sqrt{1+p^2(1)}}=\frac{1}{\cosh 2\theta_1}.$$ Then
$$C=-\frac{1}{k}\arctanh\frac{1}{\cosh 2\theta_0}.$$ Applying the
boundary condition at $s=1$ we find the value of
$\int_{0}^{1}q(\tau)\,d\tau$ as
$$\int_{0}^{1}q(\tau)\,d\tau=-\frac{1}{k}\Big(\arctanh\frac{1}{\cosh 2\theta_1}
+\arctanh\frac{1}{\cosh 2\theta_0}\Big).$$ Since, moreover,
$q(0)=\cot\psi_0$, $q(1)=\cot\psi_1$, Lemma~\ref{l3} implies the
existence of a smooth function $q(s)$ satisfying the above relation.

The function $p(s)$ can be defined by
$$\frac{1}{\sqrt{1+p^2(s)}}=-\tanh\Big[k\int_{0}^{s}q(\tau)\,d\tau-\arctanh\frac{1}{\cosh 2\theta_0}\Big].$$

The curve
$c(s)=\big(\varphi(s),\psi(s),\theta(s))=(\varphi_0+s(\varphi_1-\varphi_0),\arccot
q(s)),\frac{1}{2}\arcsinh p(s)\big)$ is the desired horizontal curve. 
\qed
\end{pf}
\medskip

\begin{remark}
Of course, the proof is given for a particular parametrisation  by a linear function $\varphi$.
One may easily  modify this proof for an arbitrary smooth function $\varphi$ obtaining a wider class of smooth horizontal curves.
\end{remark}

Some of the points on $AdS$ can be connected by a curve that maintain one of the coordinate constant. 
\medskip

\begin{theorem}\label{t5}
If $P=P(\varphi_0,\psi,\theta_0)$ and $Q=Q(\varphi_1,\psi,\theta_1)$ with 
\begin{equation}\label{psi}\psi =\arccot\Big(\ln\frac{\tanh\theta_1}{\tanh\theta_0}/\big(\varphi_0-\varphi_1\big)\Big)\end{equation} are two points that can be connected, then there is a smooth horizontal curve
joining $P$ and $Q$ with the constant $\psi$-coordinate given by~\eqref{psi}.  
\end{theorem}

\begin{pf}
Let $c=c(\varphi,\psi,\theta)$ be a horizontal curve with the constant $\psi$-coordinate. Then it satisfies the equation~\eqref{e1} that in this case we write as $$\cot\psi\ d\varphi=\frac{d(2\theta)}{\sinh 2\theta}.$$ Integrating yields $$\cot\psi\int_{\theta_0}^{\theta} d\varphi=\int_{\theta_0}^{\theta}\frac{d(2\theta)}{\sinh 2\theta}\quad\Rightarrow$$ \begin{equation}\label{e13}\cot\psi\big(\varphi(\theta)-\varphi(\theta_0)\big)=\ln\tanh\theta-\ln\tanh\theta_0.\end{equation} For $\theta=\theta_1$ we get formula~\eqref{psi} for the value of $\psi$. Solving~\eqref{e13} with respect to $\varphi(\theta)$ we get $$\varphi(\theta)=\varphi_0+\frac{\ln\big(\tanh\theta/\tanh\theta_0\big)}{\cot\psi}$$ with $\psi$ given by~\eqref{psi}. Finally, the horizontal curve joining the points $P$ and $Q$ satisfies the equation $$(\varphi,\psi,\theta)=\Big(\varphi_0+\frac{\ln\big(\tanh\theta/\tanh\theta_0\big)}{\cot\psi},\psi,\theta\Big).$$ 
\qed
\end{pf}

Upon solving the problem of the connectivity of two arbitrary points by a horizontal curve we are interested in determining its character:  timelikeness (spacelikeness or lightlikeness).  It is not an easy problem. We are able to present some particular examples showing its complexity. Let us start with the following remark.
\medskip
 
\begin{remark}\label{rem1}
If $P,Q\in AdS$ are two points connectable only by a family of smooth timelike (spacelike, lightlike) curves, then  smooth horizontal curves (its existence is known by the preceding theorem) joining $P$ and $Q$ are timelike (spacelike, lightlike).
\end{remark}
\medskip

Indeed, let $\Omega_{P,Q}$ be a family of smooth timelike (lightlike) curves connecting $P$ and~$Q$. 
If $\delta(s)\in\Omega_{P,Q}$, then its velocity vector
$\dot\delta(s)$ can be written in the left-invariant 
basis $T,X,Y$ as
$$\dot\delta(s)=\alpha(s)T(\delta(s))+\beta(s)X(\delta(s))+\gamma(s)Y(\delta(s))$$
with $\langle\dot\delta(s),\dot\delta(s)\rangle=-\alpha^2+\beta^2+\gamma^2< 0(=0)$. If moreover, it is horizontal, then $\gamma=0$. Therefore,
$-\alpha^2+\beta^2< 0(=0)$, and the horizontal curve connecting $P$ and $Q$ is 
timelike (lightlike).

If the points $P$ and $Q$ are  connectable only by   a family of spacelike curves, then the inequality
$-\alpha^2+\beta^2>\gamma^2$ holds for them. It implies $-\alpha^2+\beta^2>0$ for a horizontal curve. We conclude that 
in this case the horizontal curve is still  spacelike. 

Making use of~\eqref{palpha} and~\eqref{pbeta} as well as parametrisation~\eqref{e2} we calculate the square of the velocity vector for a horizontal curve in terms of the variables $\varphi$, $\psi$, $\theta$ as

\begin{equation}\label{sqn}
-\alpha^2+\beta^2= -\dot\varphi^2-\dot\psi^2+4\dot\theta^2-2\dot\varphi\dot\psi\cosh 2\theta.
\end{equation}
 
We present some particular timelike, spacelike, and lightlike solutions of~(\ref{e1}). 

\bigskip

\noindent
{\bf Example 2.} Let $\dot\varphi=0$.  Then, $\varphi\equiv \varphi_0$ is  constant. In order to satisfy~\eqref{e1} we have two options:
\begin{itemize}
\item[2.1]{$\dot\theta =0\quad\Longrightarrow\quad\theta\equiv \theta_0\quad\text{is  constant}$. Then $|\dot c|^2=-\dot\psi^2\leq 0$. We conclude that all non-constant horizontal curves $c(s)=(\varphi_0,\psi(s),\theta_0)$ are timelike. The projections of these curves onto the $(x_1,x_2)$- and $(x_3,x_4)$-planes are circles. All lightlike horizontal curves are only constant ones.} 
\item[2.2]{$\psi =\pi n$, $n\in\mathbb Z$. Then $|\dot c|^2=4\dot\theta^2\geq 0$. We conclude that all non-constant horizontal curves $c(s)=(\varphi_0, \pi n,\theta(s))$, $n\in\mathbb Z$ are spacelike. The projections of these curves onto the $(x_1,x_3)$- and $(x_2,x_4)$-planes are hyperbolas. All lightlike horizontal curves are only constant ones.} 
\end{itemize}

\bigskip

\noindent
{\bf Example 3.} Let $\dot\varphi\neq 0$. We choose $\varphi$ as a parameter. Then the square of the norm of the velocity vector is 

\begin{equation}\label{sqn3}
-\alpha^2+\beta^2= -1-\dot\psi^2+4\dot\theta^2-2\dot\psi\cosh 2\theta,
\end{equation}
where the derivatives are taken with respect to the parameter $\varphi$.
The horizontality condition becomes
\begin{equation}\label{hor3}
2\dot\theta\sin\psi=\cos\psi\sinh 2\theta.
\end{equation}
 As in the previous example we consider different cases.

\begin{itemize}

\item[3.1] {Suppose $\dot\theta=0$ and assume that $\theta=\theta_0\neq 0$. Then the horizontal curves are parametrized by $c(s)=(\varphi, \frac{\pi}{2}+\pi n,\theta_0)$, $n\in\mathbb Z$. All these curves are timelike, since $|\dot c|^2=-1$. There are no lightlike or spacelike horizontal curves.}
\item[3.2] {If $\theta_0=0$, then any curve in the $(\varphi,\psi)$-plane is horizontal and timelike since $|\dot c|^2=-(1+\dot\psi)^2$.}
\item[3.3] {Suppose that $\dot\psi=0$ and $\psi\equiv \psi_0\neq \frac{\pi k}{2}$, $k\in\mathbb Z$. Then~\eqref{sqn3} and~\eqref{hor3} are simplified to 
\begin{equation}\label{sqn33}
-\alpha^2+\beta^2=-1+4\dot\theta^2,
\end{equation}
\begin{equation}\label{hor33}
\dot\theta=K\sinh 2\theta\quad\text{with}\quad K=\frac{\cot\psi_0}{2}.
\end{equation}
Let $\theta=\theta(\varphi)$ solves  equation~\eqref{hor33}. Then the horizontal curve 
\begin{equation}\label{hc33}
c(s)=(\varphi, \psi_0,\theta(\varphi)) 
\end{equation}
is timelike when $|\theta|<\frac{1}{2}\arcsinh\frac{1}{2K}$. If $|\theta|>(=)\frac{1}{2}\arcsinh\frac{1}{2K}$, then the horizontal curve~\eqref{hc33} is spacelike (lightlike).}
\end{itemize}
Thus any two points $P(\varphi_0,\psi_0,\theta_0)$, $Q(\varphi_1,\psi_1,\theta_0)$, can be connected by a piecewise smooth timelike horizontal curve. This curve consists of straight segments with constant $\varphi$-coordinates or with coordinate $\psi=\frac{\pi}{2}+\pi n$, $n\in\mathbb Z$. In the case $\theta_0=0$, this horizontal curve can be constructed to be smooth.

\section{Sub-Lorentzian geodesics}\label{geodesic}

In  Lorentzian geometry there are no curves of minimal length
because two arbitrary points can be connected by a piecewise lightlike curve.
However, there do exist timelike curves with maximal length which are timelike
geodesics~\cite{Neill}. By this reason, we are looking for the longest
curve among all horizontal timelike ones. It will be given by extremizing the action integral
$S=\frac{1}{2}\int_{0}^{1}\big(-\alpha^2(s)+\beta^2(s)\big)\,ds$
under the non-holonomic constrain $\langle xE_2,\dot c\rangle=0$.
The extremal curve will satisfy the Euler-Lagrange
system \begin{equation}\label{els}\frac{d}{ds}\frac{\partial L}{\partial\dot
c}=\frac{\partial L}{\partial c}\end{equation} with the Lagrangian
$$L(c,\dot c)=\frac{1}{2}(-\alpha^2+\beta^2)+\lambda(s)\langle xE_2,\dot
c\rangle.$$ The function $\lambda(s)$ is the Lagrange multiplier
function and the values of $\alpha$ and $\beta$ are given
by~\eqref{alpha} and~\eqref{beta}. The Euler-Lagrange system~\eqref{els} can be written in the form 
\[
\begin{array}{rll}
-\dot\alpha x_2-\dot\beta x_3 & =&2(\alpha \dot x_2+\beta\dot x_3 -\lambda \dot x_4)-\dot\lambda x_4,\\
\dot\alpha x_1-\dot\beta x_4 & =& 2(-\alpha \dot x_1+\beta\dot x_4 +\lambda \dot x_3)+\dot\lambda x_3,\\           -\dot\alpha x_4+\dot\beta x_1 & =& 2(\alpha \dot x_4-\beta\dot x_1 -\lambda \dot x_2)-\dot\lambda x_2,\\          \dot\alpha x_3+\dot\beta x_2 & =& 2(-\alpha \dot x_3-\beta\dot x_2 +\lambda \dot x_1)+\dot\lambda x_4.
\end{array}
\] for the extremal curve $c(s)=(x_1(s),x_2(s),x_3(s),x_4(s))$.
Multiplying these equations by $x_2$, $-x_1$, $-x_4$, $x_3$, respectively and then, summing them up we obtain $$-\dot\alpha=2(-\alpha\langle\dot c,N\rangle-\beta\langle\dot c, Y\rangle-\lambda\beta)=-2\lambda\beta$$
because $\langle\dot c, Y\rangle=\langle\dot c, N\rangle=0$. 
Now, multiplying the equations by $x_3,x_4,x_1,x_2$, respectively and then, summing them up we get
$$-\dot\beta=2(\alpha\langle\dot c,Y\rangle+\beta\langle\dot c, N\rangle+\lambda\alpha)=2\lambda\alpha$$ in a similar way. The values of $\alpha$ and $\beta$ are concluded to satisfy the system 
\begin{equation}\label{system}
\begin{array}{lll}\dot\alpha(s) & = 2\lambda\beta(s), \\
\dot\beta(s) & =2\lambda\alpha(s).
\end{array}
\end{equation} 

{\bf Case $\lambda(s)=0$.} In the Riemannian geometry the Schwartz inequality allows us to define the angle $\vartheta$ between two vectors $v$ and $w$ as a unique number $0\leq\vartheta\leq\pi$, such that $$\cos \vartheta=\frac{v\cdot w}{|v||w|}.$$ There is an analogous result in  Lorentzian geometry which is formulated as follows.
\medskip

\begin{proposition}\label{angle}\cite{Neill}
Let $v$ and $w$ be timelike vectors. Then,
\begin{itemize}
\item[1.]{$|\langle v,w\rangle|\geq |v||w|$ where the equality is attained if and only if $v$ and $w$ are collinear.}
\item[2.]{If $\langle v,w\rangle<0$, there is a unique number $\vartheta\geq 0$, called the hyperbolic angle between $v$ and $w$, such that $$\langle v,w\rangle=-|v||w|\cosh\vartheta.$$}
\end{itemize}
 
\end{proposition}
\medskip

\begin{theorem}\label{t3}
The family of timelike future-directed horizontal curves contains horizontal timelike future-directed geodesics $c(s)$ with the following properties
\begin{itemize}
\item[1.]{The length $|\dot c|$ is constant along the geodesic.
} \item[2.]{The inner
products $\langle T,\dot c\rangle=\alpha $, $\langle X,\dot
c\rangle=\beta$, $\langle Y,\dot c\rangle=0$ are constant along
the geodesic.}
\item[3.]{The hyperbolic angle between the horizontal time vector field $T$ and the velocity vector $\dot c$ is constant.}
\end{itemize}
\end{theorem}

\begin{pf}
The system~\eqref{system} implies $$\dot\alpha(s)=0\qquad\dot\beta(s)=0.$$ The existence of a geodesic follows from the general theory of ordinary differential equations, employing, for example, the parametrisation given for $\alpha$, $\beta$, $\gamma$ in the preceding section. Since the horizontal coordinates $\alpha(s)$ and $\beta(s)$ are constant
along the curve $c$ we conclude that $c$ is geodesic. We denote by $\alpha$ and $\beta$
its respective horizontal coordinates.

The length of the velocity vector $\dot c$ is $|-\alpha^2+\beta^2|^{1/2}$ and it is constant along the geodesic.

The second statement is obvious. Since $c(s)$ is a future-directed geodesic, we have $\langle T,\dot c\rangle<0$, and $$\cosh(\angle T,\dot c)=-\frac{\langle T,\dot c\rangle}{|T||\dot c|}=\frac{-\alpha}{\sqrt{|-\alpha^2+\beta^2|}}\quad\text{is constant}.$$ 
\qed
\end{pf}

{\bf Case $\lambda(s)\neq 0$.}  We continue to study the extremals given by the solutions of the Euler-Lagrange equation~\eqref{els}.
\medskip

\begin{lemma}\label{lem1}
Let $c(s)$ be a timelike future-directed solution of the Euler-Lagrange system~\eqref{els} with $\lambda(s)\neq 0$. Then, 
\begin{itemize}
 \item [1.]{The length $|-\alpha^2(s)+\beta^2(s)|^{1/2}$ of the velocity vector $\dot c(s)$ is constant along the solution.}
\item [2.]
{The hyperbolic angle between the curve $c(s)$ and the integral curve of the time vector field $T$ is given by $$\vartheta=\angle(\dot c, T)=-2\Lambda(s)+\theta_0,$$ where $\Lambda$ is the primitive of $\lambda$.}
\end{itemize}\end{lemma}

\begin{pf}  Multiplying the first equation of~\eqref{system} by $\alpha$, the second one by $\beta$ and subtracting, we deduce that $\alpha\dot\alpha-\beta\dot\beta=0$. This implies that $-\alpha^2+\beta^2=\langle\dot c,\dot c\rangle$ is constant. The horizontal solution is timelike if the initial velocity vector is timelike. The first assertion is proved.

Set $r=\sqrt{|-\alpha^2+\beta^2|}$. Using the hyperbolic functions we write $$\alpha(s)=-r\cosh\theta(s),\qquad\beta(s)=r\sinh\theta(s).$$ Substituting $\alpha$ and $\beta$ in~\eqref{system}, we have $$\dot\theta(s)=-2\lambda(s).$$ Denote $\Lambda(s)=\int^{s}_{0}\lambda(s)\,ds$ and write the solution of the latter equation as $\theta=-2\Lambda(s)+\theta_0$. Thus, \begin{eqnarray}\label{sol}
\alpha(s)=-r\cosh(-2\Lambda(s)+\theta_0),\quad \beta(s)=r\sinh(-2\Lambda(s)+\theta_0). 
\end{eqnarray} In order to find the value of the constant $\theta_0$ we put $s=0$ and get $\theta_0=\arctanh\frac{\beta(0)}{\alpha(0)}$. 

Let $c(s)$ be a horizontal timelike future-directed solution of~\eqref{els}. Then $\langle \dot c,T\rangle<0$ and $$\alpha=\langle \dot c,T\rangle=-|\dot c| |T|\cosh\vartheta=-r\cosh(\angle(\dot c, T)).$$ Comparing with~\eqref{sol} finishes the proof of the theorem.
\qed
\end{pf}

There is no counterpart of Proposition~\ref{angle} for spacelike vectors. Nevertheless, we obtain the following analogue of Lemma~\ref{lem1} . 
\medskip

\begin{lemma}\label{lem2}
Let $c(s)$ be a spacelike solution of the Euler-Lagrange system~\eqref{els} with $\lambda(s)\neq 0$. Then, 
\begin{itemize}
 \item [1.]{The length of the velocity vector $\dot c(s)$ is constant along the solution;}
\item [2.]
{The horizontal coordinates are expressed by~\eqref{sol}.}
\end{itemize}\end{lemma}

As the next step, we shall study the function $\Lambda(s)$. First, let us prove some useful facts.
\medskip

\begin{proposition}\label{pr2}
Let $c(s)=(x_1(s),x_2(s),x_3(s),x_4(s))$ be a horizontal timelike (spacelike) curve. Then, 
\begin{itemize}
\item [1.]{$-\dot x_1^2(s)-\dot x_2^2(s)+\dot x_3^2(s)+\dot x_4^2(s)=-\alpha^2(s)+\beta^2(s)\label{hl}$;}
\item [2.]{$\ddot c=a(s)T+b(s)X+\omega(s)Y+w(s)N$, with $a=\dot\alpha$, $b=\dot\beta$, $\omega=0$, $w=\alpha^2-\beta^2$. }
\end{itemize}

\end{proposition}

\begin{pf}
Let us write the coordinates of $\dot c(s)$ in the basis $T,X,Y,N$ as $$\dot c(s)=\alpha(s)T+\beta(s)X+\gamma(s)Y+\delta(s)N,$$ where
\[\begin{array}{lll}
&\alpha & = \langle\dot c,T\rangle=x_2\dot x_1 - x_1\dot x_2 +x_4\dot
x_3 - x_3\dot x_4, \\
&\beta & =  \langle\dot c, X\rangle=-x_3\dot x_1 - x_4\dot x_2 +x_1\dot
x_3 + x_2\dot x_4, \\
0 =& \gamma &  = \langle\dot c, Y\rangle=x_4\dot x_1-x_3\dot x_2 +x_2 \dot x_3 -x_1\dot x_4, \\
0  = &\delta &  =  \langle\dot c, N\rangle=-x_1\dot x_1-x_2\dot x_2 + x_3\dot x_3 + x_4\dot x_4.
\end{array}\]
By the direct calculation we get $$-\alpha^2+\beta^2=-\alpha^2-\delta^2+\beta^2+\gamma^2=-\dot x_1^2-\dot x_2^2+\dot x_3^2+\dot x_4^2.$$

In order to prove the second statement of the proposition we calculate 
\begin{equation*}
\dot\alpha=x_2\ddot x_1 - x_1\ddot x_2 +x_4\ddot
x_3 - x_3\ddot x_4=\langle\ddot c, T\rangle=a,
\end{equation*}
\begin{equation*}
\dot\beta=-x_3\ddot x_1 - x_4\ddot x_2 +x_1\ddot
x_3 + x_2\ddot x_4=\langle\ddot c, X\rangle=b.
\end{equation*} Differentiating the horizontality condition~\eqref{hc1}, we find 
\begin{equation*}
0=\frac{d}{ds}\langle\dot c, Y\rangle=\frac{d}{ds}\big(x_4\dot x_1-x_3\dot x_2 +x_2\dot x_3 -x_1\dot x_4 \big)=x_4\ddot x_1-x_3\ddot x_2 +x_2\ddot x_3 -x_1\ddot x_4=\langle\ddot c, Y\rangle=\omega.
\end{equation*} Then, 
\begin{eqnarray*}
0 & =\frac{d}{ds}\langle\dot c, N\rangle=\frac{d}{ds}\big(-x_1\dot x_1- x_2\dot x_2+ x_3 \dot x_3+ x_4 \dot x_4\big)=-x_1\ddot x_1- x_2 \ddot x_2+ x_3 \ddot x_3+ x_4 \ddot x_4 \\ & +(-\dot x_1^2-\dot x_2^2+\dot x_3^2+\dot x_4^2)=\langle\ddot c, N \rangle + (-\alpha^2+\beta^2)=w-\alpha^2+\beta^2,
\end{eqnarray*} by the first statement. The proof is finished. 
\qed
\end{pf}
\medskip

\begin{theorem}\label{t4}
The Lagrange multiplier $\lambda(s)$ is constant along the horizontal timelike (spacelike, lightlike) solution of the Euler-Lagrange system~\eqref{els}.  
\end{theorem}

\begin{pf}
We consider the equivalent Lagrangian function $\widehat L(x,\dot x)$, changing the length function $-\alpha^2+\beta^2$ to $-\dot x_1^2-\dot x_2^2+\dot x_3^2+\dot x_4^2$. The solutions of the Euler-Lagrange system for both Lagrangians give the same curve. Thus, the new Lagrangian is $$\widehat L(x,\dot x)=\frac{1}{2}\big(-\dot x_1^2-\dot x_2^2+\dot x_3^2+\dot x_4^2\big)+\lambda(s)\big(\dot x_1x_4-\dot x_4 x_1-\dot x_2 x_3+\dot x_3 x_2\big).$$ The corresponding Euler-Lagrange system is 
\[
\begin{array}{lll}
-&\ddot x_1 & =  -\dot\lambda x_4 - 2\lambda\dot x_4,\\
-&\ddot x_2 & =  \dot\lambda x_3 + 2\lambda\dot x_3,\\
&\ddot x_3 & =  -\dot\lambda x_2 - 2\lambda\dot x_2,\\
&\ddot x_4 & =  -\dot\lambda x_1 + 2\lambda\dot x_1.
\end{array}
\]
We multiply the first equation by $-x_4$, the second equation by $x_3$, the third one by $x_2$, and the last one by $-x_1$, finally, sum them up. This yields $$\ddot x_1 x_4-\ddot x_2 x_3 +\ddot x_3 x_3 - \ddot x_4 x_1
=\dot\lambda (x_4^2+x_3^2-x_2^2-x_1^2)+2\lambda\big(\dot x_4 x_4 +\dot x_3 x_3-\dot x_2 x_2 - \dot x_1 x_1\big)\quad\Rightarrow$$ $$\langle \ddot c,Y\rangle=-\dot\lambda+2\lambda\langle \dot c,N\rangle\quad\Rightarrow\quad\dot\lambda =0.$$ We conclude that $\lambda$ is constant along the solution.
\qed
\end{pf}

We see from the proof of Lemma~\ref{lem1} that the function $\Lambda(s)$ is just a linear function. This leads to the following property of  horizontal timelike future-directed solutions of the Euler-Lagrange system~\eqref{els}.
\medskip

\begin{corollary}
If $c(s)$ is a horizontal timelike future-directed solution of~\eqref{els}, then the hyperbolic angle between its velocity and the time vector field $T$ increases linearly in~$s$.
\end{corollary}

\section{Hamiltonian formalism}\label{hamilton}

The sub-Laplacian, which is the sum of the squares of the horizontal vector fields plays the fundamental role in sub-Riemannian geometry. The counterpart of the sub-Laplacian in the Lorentz setting is the operator \begin{eqnarray}\label{lorop}\mathcal L=\frac{1}{2}(-T^2+X^2)  =  
\frac{1}{2}\Big(& - &\big(-x_2\partial_{x_1}+x_1\partial_{x_2}+x_4\partial_{x_3}-x_3\partial_{x_4}\big)^2
\nonumber \\
& + & 
\big(x_3\partial_{x_1}+x_4\partial_{x_2}+x_1\partial_{x_3}+x_2\partial_{x_4}\big)^2\Big).
\end{eqnarray}
 Then the Hamiltonian function corresponding to the operator~\eqref{lorop} is  
\begin{eqnarray}\label{ham}
H(x,\xi) & = & \frac{1}{2}\Big(-\big( -x_2\xi_1+x_1\xi_2+x_4\xi_3-x_3\xi_4\big)^2
 +  
\big(x_3\xi_1+x_4\xi_2+x_1\xi_3+x_2\xi_4\big)^2\Big)\nonumber
\\ & = & \frac{1}{2}\big(-\tau^2+\varsigma^2\big),
\end{eqnarray} where we use the notations $\xi_k=\partial_{x_k}$,  $\tau=-x_2\xi_1+x_1\xi_2+x_4\xi_3-x_3\xi_4$, and $\varsigma=x_3\xi_1+x_4\xi_2+x_1\xi_3+x_2\xi_4$. There are close relations between the solutions of the Euler-Lagrange equation and the solutions of the Hamiltonian system $$\dot x=\frac{\partial H}{\partial\xi},\qquad \dot\xi=-\frac{\partial H}{\partial x}.$$ The solutions of the Euler-Lagrange system~\eqref{els} coincide with the projection of the solutions of the Hamiltonian system onto the Riemannian manifold. In the sub-Riemannian case the solutions coincide, if and only if, the solution of the Euler-Lagrange system is a horizontal curve. We are interested in relations of the solutions of these two systems in our situation. The Hamiltonian system admits the form
\begin{equation}\label{hams}
\left\{
\begin{array}{ll}\dot x=\frac{\partial H}{\partial\xi}=-\tau xJ + \varsigma xE_1 , \\
\dot\xi=-\frac{\partial H}{\partial x}=-\tau \xi J - \varsigma \xi E_1.\end{array}\right.
\end{equation}
\medskip

\begin{lemma}\label{lem22}
The solution of the Hamiltonian system~\eqref{hams} is a horizontal curve and
\begin{equation}\label{coef}\tau=\alpha,\qquad \varsigma=\beta,\end{equation}
where $\alpha$ and $\beta$ are given by~\eqref{alpha} and~\eqref{beta} respectively.
\end{lemma}

\begin{pf}
Let $c(s)=\big(x_1(s),x_2(s),x_3(s),x_4(s)\big)$ be a solution of~\eqref{hams}. In order to prove its horizontality we need to show that the inner product $\langle\dot x, xE_2\rangle$ vanishes. We substitute $\dot x$ from~\eqref{hams} and get $$\langle\dot x, xE_2\rangle=-\tau \langle xJ , xE_2\rangle+\varsigma \langle xE_1 , xE_2\rangle=0$$ by~\eqref{m7}.

Using the first line in the Hamiltonian system and the definitions of horizontal coordinates~\eqref{alpha} and~\eqref{beta}, we get $$\alpha=\langle\dot x, xJ \rangle=-\tau \langle xJ , xJ \rangle+\varsigma \langle xE_1 , xJ \rangle=\tau,$$ 
$$\beta=\langle\dot x, xE_1 \rangle=-\tau \langle xJ , xE_1 \rangle+\varsigma \langle xE_1 , xE_1 \rangle=\varsigma$$ from~\eqref{m7} and~\eqref{m8}.
\qed
\end{pf}

\subsection{Geodesics with constant horizontal coordinates}

Lemma~\ref{lem22} implies the following form of the Hamiltonian system~\eqref{hams}  
\begin{equation}\label{hams1}\begin{array}{lll}
\dot x_1 & = & -\alpha(-x_2)+\beta x_3,\\
\dot x_2 & = & -\alpha x_1+\beta x_4,\\
\dot x_3 & = & -\alpha x_4+\beta x_1,\\
\dot x_4 & = & -\alpha(-x_3)+\beta x_2,
\end{array}
\end{equation}
with constant $\alpha$ and $\beta$.

\subsubsection{Timelike case} In this section we are aimed at finding geodesics corresponding to the extremals (Section~\ref{geodesic}) with constant horizontal coordinates $\alpha$  and $\beta$ giving the vanishing value to the Lagrangian multiplier $\lambda$. We give an explicit picture for the base point $(1,0,0,0)$. Left shifts transport it to any other point
of $AdS$.  Without lost of generality, let us assume that $-\alpha^2+\beta^2=-1$, $\alpha=\cosh \psi$, $\beta=\sinh \psi$, where $\psi$ is a constant.

The Hamiltonian system \eqref{hams1} written for constant $\alpha$  and $\beta$ is reduced to a second-order differential equation
\begin{equation}\label{trsol}
\ddot x_k=-x_k,\quad k=1,\dots 4.
\end{equation} The general solution is given
in the trigonometric basis as $x_k=A_k\cos s+B_k\sin s$. The initial condition $x(0)=(1,0,0,0)$ defines
the coefficients $A_k$ by $A_1=1$, $A_2=A_3=A_4=0$. Returning back to the first-order system \eqref{hams1} we
calculate the coefficients $B_k$ as $B_1=0$, $B_2=-\alpha$, $B_3=\beta$, $B_4=0$. Finally, the solution
is
\begin{equation}
x_1=\cos s,\quad x_2=-\cosh \psi \sin s,\quad x_3=\sinh \psi \sin s,\quad x_4\equiv 0.\label{surface1}
\end{equation}
These timelike geodesics are closed.  Varying $\psi$ they  sweep out the one-sheet hyperboloid $x_1^2+x_2^2-x_3^2=1$ in $\mathbb R^3$.  

Let us calculate the {\it vertical line} $\Gamma$, the line corresponding to the vanishing horizontal velocity $(\alpha,\beta)$
and with the constant value $\gamma=1$, passing the base point $(1,0,0,0)$. Its parametric representation $\Gamma=\Gamma(s)$ satisfies the system
\[
\begin{array}{rclcl}
\alpha&=& x_2\dot x_1 - x_1\dot x_2 +x_4\dot x_3 - x_3\dot x_4 &=& 0, \\
\beta&= &-x_3\dot x_1 - x_4\dot x_2 +x_1 \dot x_3 + x_2\dot x_4 &=&0, \\
\gamma&= &x_4\dot x_1-x_3\dot x_2 +x_2 \dot x_3 -x_1\dot x_4 &=&1, \\
\delta&= &x_1\dot x_1+x_2\dot x_2 - x_3\dot x_3 - x_4\dot x_4 &=&0.
\end{array}
\]
The discriminant of this system calculated with respect to the derivatives as variables is (-1), and we reduce the 
system to a simple one
$$
\dot x_1=-x_4,\quad \dot x_2=x_3,\quad \dot x_3=x_2,\quad \dot x_4=-x_1,
$$
with the initial condition $\Gamma(0)=x(0)=(1,0,0,0)$. The solution is $$\Gamma(s)=(\cosh s,0,0,-\sinh s).$$
The vertical line (hyperbola) $\Gamma$ meets the surface~(\ref{surface1}) at the point (1,0,0,0) orthogonally with respect to the scalar product in $\mathbb R^{2,2}$.
 Comparing this picture with the classical sub-Riemannian case of the
Heisenberg group, we observe that in the Heisenberg case all straight line geodesics lie on the horizontal plane
$\mathbb R^2$ and the center is the third vertical axis. In our case the surface~(\ref{surface1}) corresponds to the horizontal plane, timelike geodesics correspond to the straight line Heisenberg geodesics, and $\Gamma$ corresponds
to the vertical center. 

\subsubsection{Spacelike/lightlike case} Again we consider constant horizontal coordinates $\alpha$  and $\beta$, and let us assume that $-\alpha^2+\beta^2=1$, $\alpha=\sinh \psi$, $\beta=\cosh \psi$, where $\psi$ is a constant.

The Hamiltonian system \eqref{hams1} is reduced to the second-order differential equation 
\begin{equation}\label{trsol1}
\ddot x_k=x_k,\quad k=1,\dots 4.
\end{equation}  Arguing as in the previous case we deduce the solution
passing the point (1,0,0,0) as
\begin{equation}\label{surface3}
x_1=\cosh s,\quad x_2=-\sinh \psi \sinh s,\quad x_3=\cosh \psi \sinh s,\quad x_4\equiv 0.
\end{equation}
These non-closed spacelike geodesics sweep the same hyperboloid of one sheet in $\mathbb R^3$. The vertical line $\Gamma$ meets orthogonally  each spacelike geodesic on this hyperboloid at the point (1,0,0,0).

In the lightlike case $\alpha^2=\beta^2=1$ the Hamiltonian system \eqref{hams1} has a linear solution given by
$$
x_1\equiv 1,\quad x_2=-\alpha s,\quad x_3=\beta s,\quad x_4\equiv 0,
$$
which are two straight lines on the hyperboloid, and again $\Gamma$ meets them orthogonally at the unique point (1,0,0,0).

\subsection{Geodesics with non-constant horizontal coordinates.}

If the horizontal coordinates are not constant, then we must solve the Hamiltonian system generated by the
Hamiltonian \eqref{ham}.

Fix the initial point $x^{(0)}=(1,0,0,0)$. We shall give two approaches to solve this Hamiltonian system based on a  solution in Cartesian coordinates and on a parametrization of $AdS$.
\medskip

\noindent
{\bf Solution in the Cartesian coordinates.} It is convenient to introduce auxiliary phase functions 
\[
u_1=x_1+x_2,\quad u_2=x_1-x_2,\quad u_3=x_3+x_4,\quad u_4=x_3-x_4,
\] 
and momenta
\[
\psi_1=\xi_1+\xi_2,\quad \psi_2=\xi_1-\xi_2,\quad \psi_3=\xi_3+\xi_4,\quad \psi_4=\xi_3-\xi_4.
\] 
Then the Hamiltonian \eqref{ham} admits the form
$H=(-u_4\psi_2+u_1\psi_3)(u_3\psi_1-u_2\psi_4)$, and yields the Hamiltonian system 
\begin{equation} \label{hams22}
\begin{array}{lllll}
\dot u_1 & = &u_3(-u_4\psi_2+u_1\psi_3),\quad  & &u_1(0) =  1,\\
\dot u_2 & = &-u_4(u_3\psi_1-u_2\psi_4),\quad & &u_2(0) =  1,\\
\dot u_3 & = &u_1(u_3\psi_1-u_2\psi_4),\quad & &u_3(0)  = 0,\\
\dot u_4 & = &-u_2(-u_4\psi_2+u_1\psi_3),\quad & &u_4(0)  =  0,
\end{array}
\end{equation}
for positions and
\begin{equation} \label{hams23}
\begin{array}{lllll}
\dot \psi_1 & = &-\psi_3(u_3\psi_1-u_2\psi_4),\quad & &\psi_1(0)=A,\\
\dot \psi_2 & = &\psi_4(-u_4\psi_2+u_1\psi_3),\quad & & \psi_2(0)=B,\\
\dot \psi_3 & = &-\psi_1(-u_4\psi_2+u_1\psi_3),\quad & &\psi_3(0)=C,\\
\dot \psi_4 & = &\psi_2(u_3\psi_1-u_2\psi_4),\quad & &\psi_4(0)=D,
\end{array}
\end{equation}
 for momenta with some real constants $A,B,C,$ and $D$. For $\tau$ and $ \varsigma$ constant we get simple solutions mentioned in the previous section.
We see that the system (\ref{hams22}--\ref{hams23}) has the first integrals
\[\begin{array}{ccc}
u_1 \psi_1+u_3\psi_3 & = &A,\\
u_2 \psi_2+u_4\psi_4 & = &B,\\
u_2 \psi_3-u_4\psi_1 & = &C,\\
u_1 \psi_4-u_3\psi_2 & = &D,
\end{array}
\]
and in addition, we normalize $\psi(0)$ so that the trajectories belong to $AdS$: $u_1u_2+u_3u_4=1$, and 
the Hamiltonian $H=-1$ in the timelike case, in particular, the latter implies $CD=1$. Then we can deduce the momenta as
\[
\begin{array}{ccc}
\psi_1 & = &Au_2-Cu_3,\\
\psi_2 & = &Bu_1-Du_4,\\
\psi_3 & = &Cu_1+Au_4,\\
\psi_4 & = &Du_2+Bu_3.
\end{array}
\]
Let us set the functions $p=u_4/u_1$ and $q=u_3/u_2$. Then substituting function $\psi$ in \eqref{hams22}, we get
\[
\begin{array}{ccc}
\dot p & = &-(Dp^2+(A-B)p+1/D),\quad p(0)=0,\\
\dot q & = &-(Cq^2-(A-B)q+1/C), \quad q(0)=0.
\end{array}
\]
The cases of the discriminant give the following options. 
Solving these equations for $|A-B|> 2$, we obtain 
\[
p(s)=\frac{2}{D}\frac{1-e^{-s\sqrt{(B-A)^2-4}}}{(B-A-\sqrt{(B-A)^2-4})-(B-A+\sqrt{(B-A)^2-4})e^{-s\sqrt{(B-A)^2-4}}},
\] 
\[
q(s)=\frac{2D(1-e^{-s\sqrt{(A-B)^2-4}})}{(A-B-\sqrt{(A-B)^2-4})-(A-B+\sqrt{(A-B)^2-4})e^{-s\sqrt{(A-B)^2-4}}}.
\] 
Next we use the relation $\dot u_1=-\frac{u_3}{u_2}\dot u_4$. Then, $\dot u_1(pq+1)=-\dot p q u_1$, and finally,
\[
u_1(s)=\exp \int_0^s\frac{-\dot p(t) q(t)}{p(t)q(t)+1}dt,
\]
\[
u_4(s)=p(s)\exp \int_0^s\frac{-\dot p(t) q(t)}{p(t)q(t)+1}dt.
\]
Taking into account $\dot u_2=-\dot u_3 p$, we get
\[
u_2(s)=\exp \int_0^s\frac{-\dot q(t) p(t)}{p(t)q(t)+1}dt,
\]
\[
u_3(s)=q(s)\exp \int_0^s\frac{-\dot q(t) p(t)}{p(t)q(t)+1}dt.
\]
For $A-B=2$ we get
\[
u_1=(1+s)e^{-s},\quad u_2=(1-s)e^s, \quad u_3=-Dse^s,\quad u_4=-\frac{s}{D}e^{-s},
\]
or in the original coordinates
\[
x_1=\cosh s-s\,\sinh s,\quad x_2=-\sinh s+s\,\cosh s,
\]
\[
x_3=-\frac{s}{2}\left(D e^s+\frac{e^{-s}}{D}\right),\quad x_4=-\frac{s}{2}\left(D e^s-\frac{e^{-s}}{D}\right).
\]

For $A-B=-2$ and for $|A-B|< 2$  one obtains the solution analogously in the timelike case $CD=1$. Thus we get a two-parameter $D$ and $A-B$ family of geodesics passing through the point $(1,0,0,0)$. The parameters $D$ and $A-B$ have a clear dynamical meaning. Namely, 
\[
D=-\dot u_3(0)=-(\dot x_3(0)+\dot x_4(0)),\quad C=\frac{1}{D}=-\dot u_4(0)=-(\dot x_3(0)-\dot x_4(0)),
\]
and
\[
A-B=\frac{\ddot u_3(0)}{\dot u_3(0)}=-\frac{\ddot u_4(0)}{\dot u_4(0)}=\frac{\ddot x_3(0)+\ddot x_4(0)}{\dot x_3(0)+\dot x_4(0)}=-\frac{\ddot x_3(0)-\ddot x_4(0)}{\dot x_3(0)-\dot x_4(0)}.
\]
 The spacelike case
$CD=-1$ is treated in a similar way, but we omit awkward formulas. 

\medskip

\noindent
{\bf Parametric solution.} 
We present the parametric form of timelike and spacelike geodesics starting from the point $(1,0,0,0)$. The forms of solutions with constant velocity coordinates~\eqref{surface1} and~\eqref{surface3} give us an idea of a suitable parametrization for  geodesics with different causality. 

{\it Timelike geodesics.}   

We use the parametrization in a neighborhood of $(1,0,0,0)$, given by
\begin{equation}\label{e211} 
\begin{array}{lll}
x_1 & = & \cos\phi\cosh\chi_1,  \\ 
x_2 & = & \sin\phi\cosh\chi_2, \\ 
x_3 & = & \sin\phi\sinh\chi_2,\\
x_4 & = & \cos\phi\sinh\chi_1,
\end{array}
\end{equation} 
where $\phi\in(-\frac{\pi}{2},\frac{\pi}{2})$, $\chi_1,\chi_2\in(\infty,\infty)$. We note that the timelike solution with constant velocity coordinates~\eqref{surface1} followed from this parametrization if we set $\phi=-s$, $\chi_1=0$, and $\chi_2=-\psi$. The vertical line $\Gamma$ is obtained by setting $\phi=0$, $\chi_1=-s$, and $\chi_2=0$.

In this parametrization the vector fields $T$, $X$, and $Y$ admit the form
\begin{eqnarray*} 
T & = & 2\cosh(\chi_1-\chi_2)\partial_{\phi}+\partial_{\chi_1}\tan\phi\sinh(\chi_1-\chi_2)+\partial_{\chi_2}\cotan\phi\sinh(\chi_1-\chi_2), \nonumber  \\ 
X & = & 2\sinh(\chi_1-\chi_2)\partial_{\phi}+\partial_{\chi_1}\tan\phi\cosh(\chi_1-\chi_2)+\partial_{\chi_2}\cotan\phi\cosh(\chi_1-\chi_2), \nonumber \\
Y & = & \partial_{\chi_1}-\partial_{\chi_2}.
\end{eqnarray*} The vertical direction is given by the constant vector field $Y$. Let $c(s)=(\phi(s),\chi(s),\chi_2(s))$ be a curve starting at $c(0)=(0,0,\chi_2(0))$. The horizontal coordinates~\eqref{alpha} and~\eqref{beta} with respect to given parametrization are 
\begin{eqnarray*}
\alpha & = & -\dot\phi\cosh(\chi_1-\chi_2)+\frac{1}{2}(\dot\chi_1+\dot\chi_2)\sin(2\phi)\sinh(\chi_1-\chi_2),\\
\beta & = & \dot\phi\sinh(\chi_1-\chi_2)+\frac{1}{2}(\dot\chi_1+\dot\chi_2)\sin(2\phi)\cosh(\chi_1-\chi_2).
\end{eqnarray*} Then, the square of the velocity vector $\dot c(s)$ is $$-\alpha^2+\beta^2=-\dot\phi^2+\frac{1}{4}(\dot\chi_1+\dot\chi_2)^2\sin^2(2\phi).$$ The speed is preserved along the geodesics and is equal to the initial value at the point $(1,0,0,0)$, or in our parametrization $(0,0,\chi_2(0)$. Therefore, $$\langle\dot c(0),\dot c(0) \rangle=(-\alpha^2+\beta^2)(0)=-\dot\phi^2(0),$$ and we obtain timelike geodesics starting from $(0,0,\chi_2(0))$ if $\dot\phi(0)\neq 0$, and lightlike geodesics in the limiting case $\dot\phi(0)= 0$. 

The Hamiltonian $H$ associated with the operator $$\mathcal L=\frac{1}{2}(-T^2+X^2)=\frac{1}{2}(-4\partial^2_{\phi}+\tan^2\phi\partial^2_{\chi_1}+\cotan^2\phi\partial^2_{\chi_2}+2\partial_{\chi_1}\partial_{\chi_2}),$$ becomes  $$H(\phi,\chi_1,\chi_2,\psi,\xi_1,\xi_2)=\frac{1}{2}(-4\psi^2+\xi_1^2\tan^2\phi+\xi_2^2\cotan^2\phi+2\xi_1\xi_2),$$ where we set $\partial_{\phi}=\psi$, $\partial_{\chi_1}=\xi_1$, and $\partial_{\chi_2}=\xi_2$.

The Hamiltonian system 
\begin{equation}\label{hamsyst}
\begin{array}{ccl}
\dot\chi_1 & = & \xi_1\tan^2\phi+\xi_2,  \\
\dot\chi_2 & = & \xi_2\cotan^2\phi+\xi_1,\\
\dot\phi & = & -4\psi, \\
\dot\xi_1 & = & 0, \\
\dot\xi_2 & = &  0, \\
\dot\psi & = & -\xi_1^2\frac{\tan\phi}{\cos^2\phi}+\xi_2^2\frac{\cotan\phi}{\sin^2\phi}.
\end{array}
\end{equation}
shows that $\xi_1$ and $\xi_2$ are constants. If both constants vanish, then we get $$\dot\chi_1=0,\quad \dot\chi_2=0,\quad\dot\phi = -4\psi,\quad \dot\psi=0,$$ which leads to the trivial solution~\eqref{surface1}. Since we are looking for a solution in a neighborhood of $(0,0,\chi_2(0))$, we put $\xi_2=0$. Let us solve the Hamiltonian system~\eqref{hamsyst} with the initial conditions $$\phi(0)=0,\quad\chi_1(0)=0, \quad\chi_2(0)=\chi_2^{(0)},\qquad \psi(0)=\psi^{(0)},\quad \xi_1(0)=\xi_1, \quad \xi_2(0)=0.$$
From the third and from the last equations we get $\ddot\phi=-4\dot\psi=4\xi_1^2\frac{\tan\phi}{\cos^2\phi}$. Multiplying by $\dot\phi$ and integrating we obtain \begin{equation}\label{dotphi}\dot\phi^2(s)=C^2+4\xi_1^2\tan^2\phi(s),\qquad C=\dot\phi^2(0)=16\psi^2(0).\end{equation}

Let us assume $C^2>0$. Simplifying~\eqref{dotphi}, we come to $$\frac{\cos\phi\,d\phi}{\sqrt{C^2+(4\xi_1^2-C^2)\sin^2\phi}}=\pm ds.$$ According to the sign of $4\xi_1^2-C^2$, one gets three different types of solutions.

\noindent {\it Case 1: $4\xi_1^2-C^2=0$}. Integrating from $0$ to some value of $s$ we get the solution in the form  $\sin\phi(s)=\pm |C|\,s$.

\noindent {\it Case 2: $4\xi_1^2-C^2>0$}. The solution  follows as $$\sqrt{\frac{4\xi_1^2-C^2}{C^2}}\sin\phi=\pm\sinh( s\sqrt{4\xi_1^2-C^2}).$$

\noindent {\it Case 3: $4\xi_1^2-C^2<0$}. The solution is obtained as  $$\sqrt{\frac{C^2-4\xi_1^2}{C^2}}\sin\phi=\pm \sin(s \sqrt{C^2-4\xi_1^2}).$$

In order to calculate the value of $\chi_1$, we express $\tan^2\phi$ from the Cases 1-3 and integrate the first equation of the Hamiltonian system. Observe that $\dot\chi_2=\xi_1$ is constant and $\dot\phi(0)=-4\psi^{(0)}\neq 0$. The following theorem is proved.
\medskip

\begin{theorem}
The timelike geodesics starting from the point $\phi(0)=0$, $\chi_1(0)=0$, $\chi_2(0)=\chi_2^{(0)}$ with some $\dot\phi(0)$, a constant value of $\dot\chi_2$, and an arbitrary $\dot\chi_1(s)$ satisfy the following equations:

\smallskip

\noindent if $4\dot\chi_2^2=\dot\phi^2(0)$ then
\begin{itemize}
\item[$\bullet$]{$\sin\phi(s)=\pm |C|s$,}
 \item[$\bullet$]{$\chi_1(s)=-\dot\chi_2 s+\frac{\dot\chi_2}{2\dot\phi(0)}\ln\Big|\frac{1+\dot\phi(0) s}{1-\dot\phi(0) s}\Big|$,}
\item[$\bullet$] {$\chi_2(s)=\dot\chi_2 s+\chi_2^{(0)}$;}
\end{itemize}

\smallskip

\noindent if $4\dot\chi_2^2>\dot\phi^2(0)$ then
\begin{itemize}
\item[$\bullet$]{$\sin\phi(s)=\pm \sqrt{\frac{\dot\phi^2(0)}{4\dot\chi_2^2-\dot\phi^2(0)}}\sinh\Big(s\sqrt{4\dot\chi_2^2-\dot\phi^2(0)}\Big)$,}
 \item[$\bullet$]{$\chi_1(s)=-\dot\chi_2 s+\dot\chi_2\int_0^{s}\frac{4\dot\chi_2^2-\dot\phi^2(0)}{4\dot\chi_2^2-\dot\phi^2(0)\cosh^2(s\sqrt{4\dot\chi_2^2-\dot\phi^2(0)})}ds$,}
\item[$\bullet$] {$\chi_2(s)=\dot\chi_2 s+\chi_2^{(0)}$;}
\end{itemize}

\smallskip

\noindent and if f $4\dot\chi_2^2<\dot\phi^2(0)$ then
\begin{itemize}
\item[$\bullet$]{$\sin\phi(s)=\pm \sqrt{\frac{\dot\phi^2(0)}{\dot\phi^2(0)-4\dot\chi_2^2}}\sin\Big(s \sqrt{\dot\phi^2(0)-4\dot\chi_2^2}\Big)$,}
 \item[$\bullet$]{$\chi_1(s)=-\dot\chi_2 s+\dot\chi_2\int_0^s\frac{\dot\phi^2(0)-4\dot\chi_2^2}{\dot\phi^2(0)\cos^2(s\sqrt{\dot\phi^2(0)-4\dot\chi_2^2})-4\dot\chi_2^2}ds$,}
\item[$\bullet$] {$\chi_2(s)=\dot\chi_2 s+\chi_2^{(0)}$.}
\end{itemize}
\end{theorem}
\medskip
The integrals can be easily calculated and they involve  trigonometric and hyperbolic functions, and depend on the relations between $4\dot\chi_2^2$, $\dot\phi^2(0)$. 

\medskip

{\it Spacelike geodesics.}

We use another parametrizaition in a neighborhood of $(1,0,0,0)$ suitable in this case
\begin{equation}\label{e212} 
\begin{array}{rcl}
x_1 & = & \cosh\phi\cosh\chi_1,   \\ 
x_2 & = & \sinh\phi\cosh\chi_2,  \\ 
x_3 & = & \sinh\phi\sinh\chi_2,\\
x_4 & = & \cosh\phi\sinh\chi_1,
\end{array}
\end{equation} where $\phi,\chi_1,\chi_2\in(-\infty,\infty)$. Observe that the spacelike solution with constant velocity coordinates~\eqref{surface3} follows from this parametrization if we set $\phi=s$, $\chi_1=0$ and $\chi_2=-\psi$. The vertical line $\Gamma$ is obtained as previously, by setting $\phi=0$, $\chi_1=-s$, and $\chi_2=0$.

The vector fields $T$, $X$, and $Y$ become 
\begin{eqnarray*} 
T & = & 2\sinh(\chi_1-\chi_2)\partial_{\phi}-\partial_{\chi_1}\tan\phi\cosh(\chi_1-\chi_2)+\partial_{\chi_2}\cotan\phi\cosh(\chi_1-\chi_2), \nonumber  \\ 
X & = & 2\cosh(\chi_1-\chi_2)\partial_{\phi}-\partial_{\chi_1}\tan\phi\sinh(\chi_1-\chi_2)+\partial_{\chi_2}\cotan\phi\sinh(\chi_1-\chi_2), \nonumber \\
Y & = & \partial_{\chi_1}-\partial_{\chi_2}.
\end{eqnarray*} The vertical direction is again given by a constant vector field $Y$. Let $c(s)=(\phi(s),\chi(s),\chi_2(s))$ be a curve such that $c(0)=(0,0,\chi_2(0))$. The horizontal coordinates~\eqref{alpha} and~\eqref{beta} with respect to this parametrizaition are
\[
\begin{array}{ccl}
\alpha & = & \dot\phi\sinh(\chi_1-\chi_2)-\frac{1}{2}(\dot\chi_1+\dot\chi_2)\sinh(2\phi)\cosh(\chi_1-\chi_2),\\
\beta & = & \dot\phi\cosh(\chi_1-\chi_2)-\frac{1}{2}(\dot\chi_1+\dot\chi_2)\sinh(2\phi)\sinh(\chi_1-\chi_2).
\end{array} 
\]
Then the square of the velocity vector $\dot c$ is $$-\alpha^2+\beta^2=\dot\phi^2-\frac{1}{4}(\dot\chi_1+\dot\chi_2)^2\sinh^2(2\phi).$$ Since the speed is preserved along the geodesics, it is equal to $\dot\phi^2(0)$, and we obtain  spacelike geodesics starting from $(0,0,\chi_2(0))$ for $\dot\phi(0)\neq 0$. 

The Hamiltonian $H$ associated with the operator $$\mathcal L=\frac{1}{2}(-T^2+X^2)=\frac{1}{2}(4\partial^2_{\phi}-\tanh^2\phi\partial^2_{\chi_1}-\cotanh^2\phi\partial^2_{\chi_2}+2\partial_{\chi_1}\partial_{\chi_2})$$ becomes $$H(\phi,\chi_1,\chi_2,\psi,\xi_1,\xi_2)=\frac{1}{2}(4\psi^2-\xi_1^2\tan^2\phi-\xi_2^2\cotan^2\phi+2\xi_1\xi_2),$$ where we set $\partial_{\phi}=\psi$, $\partial_{\chi_1}=\xi_1$, and $\partial_{\chi_2}=\xi_2$.

As in the previous case, the Hamiltonian system 
\begin{equation}\label{hamsyst2}
\begin{array}{ccl}
\dot\chi_1 & = & -\xi_1\tanh^2\phi+\xi_2 , \\
\dot\chi_2 & = & -\xi_2\cotanh^2\phi+\xi_1,\\
\dot\phi & = & 4\psi, \\
\dot\xi_1 & = & 0, \\
\dot\xi_2 & = &  0,\\
\dot\psi & = & \xi_1^2\frac{\tanh\phi}{\cosh^2\phi}-\xi_2^2\frac{\cotanh\phi}{\sinh^2\phi}.
\end{array}
\end{equation}
gives that $\xi_1$ and $\xi_2$ are constants. If both constants vanish, we get $$\dot\chi_1=0,\quad \dot\chi_2=0,\quad\dot\phi = -4\psi,\quad \dot\psi=0$$ which leads to the spacelike trivial solution. Setting $\xi_2=0$, we solve the Hamiltonian system~\eqref{hamsyst2} with the initial conditions $$\phi(0)=0,\quad\chi_1(0)=0, \quad\chi_2(0)=\chi_2^{(0)},\qquad \psi(0)=\psi^{(0)},\quad \xi_1(0)=\xi_1, \quad \xi_2(0)=0.$$
An analogue of~\eqref{dotphi} is \begin{equation}\label{dotphi2}\dot\phi^2(s)=C^2+4\xi_1^2\tanh^2\phi(s),\qquad C=\dot\phi^2(0)=16\psi^2(0)\neq 0.\end{equation} Arguing as in the timelike case, we prove the following statement. 
\medskip

\begin{theorem}
The spacelike geodesics starting from the point $\phi(0)=0$, $\chi_1(0)=0$, $\chi_2(0)=\chi_2^{(0)}$ with some $\dot\phi(0)$, a constant value of $\dot\chi_2$, and an arbitrary $\dot\chi_1(s)$ have the following equations:
\begin{eqnarray*}
\sinh\phi(s) & = & \pm \sqrt{\frac{\dot\phi^2(0)}{\dot\phi^2(0)+4\dot\chi_2^2}}\sinh(s \sqrt{\dot\phi^2(0)+4\dot\chi_2^2}),\\c
\chi_1(s) & = & -\dot\chi_2 s+\frac{\dot\chi_2}{2|\dot\chi_2|}\arccotanh\Big(\sqrt{\frac{\dot\phi^2(0)+4\dot\chi_2^2}{4\dot\chi_2^2}}\cotan\big(s\sqrt{\dot\phi^2(0)+4\dot\chi_2^2}\big)\Big),\\
\chi_2(s) & = & \dot\chi_2 s+\chi_2^{(0)}.
\end{eqnarray*}
\end{theorem}

\section{Geodesics with respect to the distribution $\mathcal D=\spn\{X,Y\}$}\label{xy}

This case reveals the sub-Riemannian nature of such a distribution. In principle, one can easily modify the classical results from sub-Riemannian geometry (Chow-Rashevskii theorem, in particular). However we prefer to modify our own results proved in previous sections to show some particular features and to compare with the sub-Lorentzian case defined by the distribution $\mathcal D=\spn\{T,X\}$. 
\medskip

\begin{lemma}\label{l11} A curve
$c(s)=(x_1(s),x_2(s),x_3(s),x_4(s))$ is horizontal with respect to
the distribution $\mathcal D=\spn\{X,Y\}$, if and only if,
\begin{equation}\label{horcon}x_2\dot x_1-x_1\dot x_2+x_4\dot x_3-x_3\dot x_4=0 \qquad \text{or}\qquad \langle x J ,\dot c\rangle=0.\end{equation}
\end{lemma}

\begin{pf}
The tangent vector to a curve $c(s)=(x_1(s),x_2(s),x_3(s),x_4(s))$
written in the left-invariant basis is of the form $$\dot
c(s)=\alpha T+\beta X + \gamma Y.$$ Then
$$\alpha=\langle\dot c, T\rangle=\mathcal I\dot c\cdot T=\dot x_1 x_2-\dot x_2 x_1+\dot x_3 x_4-\dot x_4
x_3.$$ We conclude that $\alpha=0$, if and only
if, ~\eqref{horcon} holds.
\qed
\end{pf}

In this case a curve is horizontal, if and only if, its velocity
vector is orthogonal to the vector field $T$. The left-invariant
coordinates $\beta (s)$ and $\gamma(s)$ of a horizontal curve
$c(s)=(x_1(s),x_2(s),x_3(s),x_4(s))$ are
\begin{equation}\label{beta1}
\beta=\langle\dot c, X\rangle=-x_3\dot x_1 - x_4\dot x_2 +x_1\dot
x_3 + x_2\dot x_4=\langle xE_1 ,\dot c\rangle.
\end{equation}
\begin{equation}\label{gamma1}
\gamma=\langle\dot c, Y\rangle=-x_4\dot x_1 + x_3\dot x_2 - x_2\dot
x_3 + x_1\dot x_4=\langle xE_2,\dot c\rangle.
\end{equation}

The form $w=-x_2 dx_1 + x_1 dx_2 - x_4 dx_3 + x_3 dx_4=-\langle
xJ ,dx \rangle$ is a contact form for the horizontal distribution
$\mathcal D=\spn\{X,Y\}$. Indeed $$w(N)=0,\quad w(T)=1,\quad
w(X)=0,\quad w(Y)=0.$$ Thus, $\ker w=\spn \{N,X,Y\}$, The
horizontal distribution can be defined as follows $$\mathcal
D=\{V\in TAdS:\ w(V)=0\},\quad\text{or}\quad \mathcal D=\ker w\cap
TAdS.$$

The length $l(c)$ of a horizontal curve $c(s):\ [0,1]\to AdS$ is
given by
$$l(c)=\int_{0}^{1}\langle\dot c(s),\dot c(s)\rangle^{1/2}\,ds
=\int_{0}^{1}\big(\beta^2(s)+\gamma^2(s)\big)^{1/2}\,ds.$$ The
restriction of  the non-degenerate metric
$\langle\cdot ,\cdot\rangle$ onto the horizontal distribution
$\mathcal D\subset TAdS$ gives a positive-definite metric that we still denote by $\langle\cdot ,\cdot\rangle_{\mathcal D}$. Thus from now on, we shall work only with one type of the curves (that we shall call simply horizontal curves), since the horizontality condition requires the vanishing coordinate function of the vector field $T$. 

\subsection{Existence of horizontal curves}

The following theorem is an analogue to Theorem~\ref{t2} proved for the distribution $\mathcal D=\spn\{T,X\}$ in Section 4.

\medskip

\begin{theorem}
Let $P$, $Q\in AdS$ be arbitrary given points. Then there is a smooth  horizontal curve connecting $P$ with $Q$.
\end{theorem}

\begin{pf}
We use parametrisation~\eqref{e2}, in which  the horizontality condition for a curve $c(s)$ is expressed by~\eqref{palpha} as $$\dot\psi+\dot\varphi\cosh 2\theta=0.$$  This equation is to be sold for the initial conditions $$c(0)=P,\quad\text{or}\quad\varphi(0)=\varphi_0,\quad\psi(0)=\psi_0,\quad\theta(0)=\theta_0,$$ $$c(1)=Q,\quad\text{or}\quad\varphi(1)=\varphi_1,\quad\psi(1)=\psi_1,\quad\theta(1)=\theta_1.$$  Let $\psi=\psi(s)$ be a smooth arbitrary function with $\dot\psi(0)=\lim\limits_{s\to 0^+}\dot\psi(s)$ and $\dot\psi(1)=\lim\limits_{s\to 1^-}\dot\psi(s)$. Set $2\theta(s)=\arccosh p(s)$. Then the equation~\eqref{palpha} admits the form $$\dot\varphi=-\frac{\dot\psi}{\cosh 2\theta}=-\frac{\dot\psi}{p(s)}\quad\Rightarrow\quad\varphi(s)=-\int_{0}^{s}\frac{\dot\psi(s)\,ds}{p(s)}+\varphi(0).$$ Denote $q(s)=\frac{\dot\psi(s)}{p(s)}$. Since $q(0)=\frac{\dot\psi(0)}{\cosh 2\theta_0}$, $q(1)=\frac{\dot\psi(1)}{\cosh 2\theta_1}$, and $\int_{0}^{1}q(s)\,ds=\varphi_0-\varphi_1$ applying Lemma~\ref{l3} we conclude that there exists such a smooth function $q(s)$. The function $p(s)$ is found as $p(s)=\frac{\dot\psi(s)}{q(s)}$. We get a curve $c(s)=(\varphi(s),\psi(s),\theta(s))$ with
\[
\begin{array}{lll}
\psi& = &\psi(s),\\
\varphi(s)& =& -\int_{0}^{s}\frac{\dot\psi(s)\,ds}{p(s)}+\varphi(0),\\
\theta(s)& = &\frac{1}{2}\arccosh p(s).
\end{array}
\]
\qed
\end{pf}
\medskip

\begin{remark} Observe that in the general Chow-Rashevskii theorem smoothness was not concluded.\end{remark}

\medskip

\begin{theorem}
Given two arbitrary points $P=P(\varphi_0,\psi_0,\theta_0)$ and $Q=Q(\varphi_1,\psi_1,\theta_0)$ with $2\theta_0=\arccosh\frac{\psi_1-\psi_0}{\varphi_0-\varphi_1}$, there is a horizontal curve with the constant $\theta$-coordinate connecting $P$ with~$Q$. 
\end{theorem}

\begin{pf}
If the $\theta$-coordinate is constant, then the governing equation is $$\dot\psi=-\dot\varphi\cosh 2\theta_0\quad\Rightarrow\quad\psi(s)=-\varphi(s)\cosh 2\theta_0+C.$$
Applying the initial conditions $$c(0)=\big(\varphi_0,\psi_0,\theta_0\big),\quad\text{and}\quad c(1)=\big(\varphi_1,\psi_1,\theta_0\big),$$ we find $$2\theta_0=\arccosh\Big(\frac{\psi_1-\psi_0}{\varphi_0-\varphi_1}\Big),\qquad C=\psi_0+\varphi_0\frac{\psi_1-\psi_0}{\varphi_0-\varphi_1}.$$ Therefore, for any parameter $\varphi$, the horizontal curve $$c(s)=\Big(\varphi,\psi_0+\big(\varphi(0)-\varphi\big)\frac{\psi_1-\psi_0}{\varphi_0-\varphi_1},\theta_0\Big),\qquad 2\theta_0=\arccosh\frac{\psi_1-\psi_0}{\varphi_0-\varphi_1},$$ joins the points $P=P(\varphi_0,\psi_0,\theta_0)$ and $Q=Q(\varphi_1,\psi_1,\theta_0)$.
\qed
\end{pf}

\subsection{Lagrangian formalism}

Dealing with $\mathcal D=\spn\{X,Y\}$ and a positive-definite metric $\langle\cdot ,\cdot\rangle_{\mathcal D}$ on it, one might compare with the geometry generated by the sub-Riemannian distribution on sphere $S^3$ in~\cite{CChM}. The minimizing length curve can be found by minimizing the action integral $$S=\frac{1}{2}\int_{0}^{1}(\beta^2(s)+\gamma^2(s))\,ds$$ under the non-holonomic constrain $\alpha=\langle \dot c,xJ \rangle=0$. The corresponding Lagrangian is \begin{equation}\label{lagr2}L(c,\dot c)=\frac{1}{2}\big(\beta^2(s)+\gamma^2(s)\big)+\lambda(s)\alpha(s).\end{equation} The extremal curve is given by the solution of the Euler-Lagrange system~\eqref{els} with the Lagrangian~\eqref{lagr2}.

Let us make some preparatory calculations. Write the system~\eqref{els} for the Lagrangian~\eqref{lagr2} as the follows
\[
\begin{array}{rcl}
2\beta\dot x_3+2\gamma\dot x_4-2\lambda\dot x_2+\dot\beta x_3+\dot\gamma x_4-\dot\lambda x_2 & = & 0,\\
2\beta\dot x_4-2\gamma\dot x_3+2\lambda\dot x_1+\dot\beta x_4-\dot\gamma x_3+\dot\lambda x_1 & = & 0,\\
-2\beta\dot x_1+2\gamma\dot x_2-2\lambda\dot x_4-\dot\beta x_1+\dot\gamma x_2-\dot\lambda x_4 & = & 0,\\
-2\beta\dot x_2-2\gamma\dot x_1+2\lambda\dot x_3-\dot\beta x_2-\dot\gamma x_1+\dot\lambda x_3 & = & 0.
\end{array}
\]
Multiply the equations by $x_3,x_4,x_1$, and $x_2$, respectively and sum them up. We get 
\[
\begin{array}{rcl}\label{bg}2\beta\langle \dot c,N\rangle-2\gamma\langle\dot c, T\rangle-2\lambda\langle\dot c, Y\rangle-\dot\beta+0\dot\gamma+0\dot\lambda=0\ \ & \Rightarrow &\ \ \dot\beta=2\lambda\gamma,\\
2\beta\langle \dot c,T\rangle-2\gamma\langle\dot c, N\rangle+2\lambda\langle\dot c, X\rangle+0\dot\beta-\dot\gamma+0\dot\lambda=0\ \ & \Rightarrow & \ \ \dot\gamma=2\lambda\beta.
\end{array}
\]
Let us consider two cases.

{\bf Case $\lambda(s)=0$.} In this case  equation~\eqref{bg} admits the form 
\begin{equation}\label{bg1}\dot\beta=0,\qquad\dot\gamma=0,\end{equation}
and we deduce the following theorem.
\medskip

\begin{theorem}
There are horizontal geodesics with the following properties:
\begin{itemize}
 \item [1.]{The coordinates $\alpha=\langle \dot c, T\rangle=0$, $\beta=\langle\dot c, X\rangle$, and $\gamma=\langle\dot c, Y\rangle$ are constant;}
\item [2.]{The length $|\dot c|$ along the geodesics;}
\item [3.]{The angles between the velocity vector and horizontal frame is constant along along the geodesic.}
\end{itemize}
\end{theorem}

\begin{pf}
Taking into account the solution of~\eqref{bg1}, we denote $\beta(s)=\beta$ and $\gamma(s)=\gamma$. Then the length of the velocity vector $|\dot c|=\sqrt{\beta^2+\gamma^2}$ is constant. 

Since $\langle\dot c, X\rangle=\langle\dot c, X\rangle_{\mathcal D}=|\dot c|_{\mathcal D}|X|_{\mathcal D}\cos(\angle\dot c,X)$, $\langle\dot c, Y\rangle=\langle\dot c, Y\rangle_{\mathcal D}=|\dot c|_{\mathcal D}|Y|_{\mathcal D}\cos(\angle\dot c,Y)$, we have $$\cos(\angle\dot c,X)=\frac{\beta}{\sqrt{\beta^2+\gamma^2}},\quad\cos(\angle\dot c,Y)=\frac{\gamma}{\sqrt{\beta^2+\gamma^2}},$$ that proves the third  assertion.
\qed
\end{pf}
\medskip

{\bf Case $\lambda(s)\neq 0$.}
\begin{theorem} There are horizontal geodesics with the following properties:
 \begin{itemize}
  \item [1.]{The velocity vector $|\dot c|$ of a geodesic is constant along the geodesic;}
  \item [2.]{The angles between the velocity vector and the horizontal frame are given by $$\angle\dot c,X=cs+\theta_0,\qquad \angle\dot c,Y=\frac{\pi}{2}-cs+\theta_0.$$}
 \end{itemize}
\end{theorem}

\begin{pf} Since
\begin{equation}\label{gb3}\dot\beta=2\lambda\gamma,\quad\dot\gamma=2\lambda\beta\end{equation} implies $\frac{d}{ds}{\big(\beta^2+\gamma^2\big)}=0$, we conclude, that the length of the velocity vector $|\dot c|$ is constant. Taking into account positivity of $\beta^2+\gamma^2$ let us denote it by $r^2$. Set $\beta=r\cos\theta(s)$ and $\gamma=r\sin\theta(s)$. Substituting them in~\eqref{gb3}, we get $$\dot\theta(s)=2\lambda(s)\ \ \Rightarrow\ \ \theta(s)=2\int\lambda(s)\,ds+\theta_0.$$ 

Let us find the function $\lambda(s)$.  Observe that $$\beta^2+\gamma^2=-\dot x_1^2-\dot x_2^2+\dot x_3^2+\dot x_4^2.$$ It can be shown similarly to the proof of Proposition~\ref{pr2}, having  $\alpha=\delta=0$. By the direct calculation (see also Proposition~\ref{pr2}) we show that $$\langle\ddot c,T\rangle=\frac{d}{ds}\langle\dot c,T\rangle=0.$$ Now, we consider an equivalent to~\eqref{lagr2} extremal problem with the Lagrangian
\begin{equation}\label{lagr3}
 \widehat L(c,\dot c)=\frac{1}{2}\big(-\dot x_1^2-\dot x_2^2+\dot x_3^2+\dot x_4^2\big)+\lambda(s)\langle\dot c, T\rangle.
\end{equation}
The Euler-Lagrange system admits the form 
\[
\begin{array}{rl}
-\ddot x_1&=-2\lambda\dot x_2-\dot\lambda x_2,\\
-\ddot x_2&=2\lambda\dot x_1+\dot\lambda x_1,\\
\ddot x_3&=-2\lambda\dot x_4-\dot\lambda x_4,\\
\ddot x_4&=2\lambda\dot x_3+\dot\lambda x_3.
\end{array}
\]
Multiplying these equations by $x_2,-x_1,-x_4,x_3$ respectively and then, summing them up, we obtain $$-\langle\ddot c, T\rangle=2\lambda\langle\dot c, N\rangle-\dot\lambda.$$ This allows us to conclude, that the function $\lambda(s)$ is constant along the solution of the Euler-Lagrange equation that yields the second assertion of the theorem. 
\qed
\end{pf}

\subsection{Hamiltonian formalism}

The sub-Laplacian is $\mathcal L=X^2+Y^2$ and the corresponding Hamiltonian function is
$$H(x,\xi)=\frac{1}{2}\Big(\big(x_3\xi_1+x_4\xi_2+x_1\xi_3+x_2\xi_4\big)^2+\big(x_4\xi_1-x_3\xi_2-x_2\xi_3+x_1\xi_4\big)^2\Big)=\frac{1}{2}(\varsigma^2+\kappa^2).$$ The Hamiltonian system is written as
\begin{equation}\label{hams3}
\begin{array}{rll}
\dot x&=\frac{\partial H}{\partial\xi}&=\varsigma xE_1 +\kappa xE_2, \\
\dot \xi&=-\frac{\partial H}{\partial x}&=-\varsigma \xi E_1  -\kappa \xi E_2,
\end{array}
\end{equation}
As in the previous section we are able to prove the following proposition.
\medskip

\begin{proposition}\label{prop}
The solution of the Hamiltonian system is a horizontal curve and $$\varsigma=\beta,\qquad \kappa=\gamma.$$
\end{proposition}
\medskip

\begin{corollary}
The Hamiltonian function is the energy $H(x,\xi)= \frac{1}{2}(\beta^2+\gamma^2)$.
\end{corollary}

\subsection{Geodesics with constant horizontal coordinates.}
 In this section we consider constant horizontal coordinates $\beta$  and $\gamma$. 
 Making use of Proposition~\ref{prop} we write the first line of the Hamiltonian system~\eqref{hams3} in the form.
\begin{equation}\label{hams4}
\begin{array}{rl}
\dot x_1&=\beta x_3+\gamma x_4,\\
\dot x_2&=\beta x_4-\gamma x_3, \\
\dot x_3&=\beta x_1-\gamma x_2, \\
\dot x_4&=\beta x_2+\gamma x_1, 
\end{array}
\end{equation}
 We give an explicit picture for the base point $(1,0,0,0)$.  Without lost of generality, let us assume that $\beta^2+\gamma^2=1$, $\beta=\cos \psi$, $\gamma=\sin \psi$, where $\psi$ is a constant.

The Hamiltonian system \eqref{hams4} written for constant $\beta$  and $\gamma$ is reduced to a second-order differential equation
\begin{equation}\label{trsol2}
\ddot x_k=x_k,\quad k=1,\dots 4.
\end{equation} The general solution is given
in the hyperbolic basis as $x_k=A_k\cosh s+B_k\sinh s$. The initial condition $x(0)=(1,0,0,0)$ defines
the coefficients $A_k$ by $A_1=1$, $A_2=A_3=A_k=0$. Returning back to the first-order system \eqref{hams4} we
calculate the coefficients $B_k$ as $B_1=0$, $B_2=0$, $B_3=\beta$, $B_4=\gamma$. Finally, the solution
is
\begin{equation}
x_1=\cosh s,\quad x_2\equiv 0,\quad x_3=\cos \psi \sinh s,\quad x_4= \sin\psi \sinh s.\label{surface2}
\end{equation}
 Varying $\psi$ they  sweep out the two-sheet hyperboloid $x_1^2-x_3^2-x_4^2=1$ in $\mathbb R^3$.  
 We use only one sheet containing the point $(1,0,0,0)$. Geodesics are hyperbolas passing this point.

The vertical line corresponds to the vanishing horizontal velocity $(\beta, \gamma)$
and with the constant value $\alpha=1$, passing the base point $(1,0,0,0)$. The solution is $$\Gamma(s)=(\cos s,\sin s, 0,0).$$
The vertical line (circle) $\Gamma$ meets the surface~(\ref{surface2}) at the point (1,0,0,0) orthogonally with respect to the scalar product in $\mathbb R^{2,2}$.

\subsection{Geodesics with non-constant horizontal coordinates.} 

If the horizontal coordinates are not constant, then we must solve the Hamiltonian system generated by the above
Hamiltonian.
\medskip

\noindent
{\bf Solution in the Cartesian coordinates.}
Fix the initial point $x^{(0)}=(1,0,0,0)$. In the Cartesian case it is convenient to introduce complex coordinates $z=x_1+ix_2$, $w=x_3+ix_4$, $\varphi=\xi_1+i\xi_2$, and $\psi=\xi_3+i\xi_4$. Hence, the Hamiltonian admits the form $H=|z\bar\psi+\bar w\varphi |^2$. The corresponding Hamiltonian system becomes
\begin{equation*}
\begin{array}{lllll}
\dot z & = &  w(z\bar{\psi}+\bar{w} \varphi),\quad & &z(0)=1, \\
\dot w & = & z(\bar{z}\psi+w\bar{\varphi}), \quad & &w(0)=0,\\
\dot {\bar{\varphi}} & = & -\bar{\psi}(\bar{z}\psi+w\bar{\varphi}), \quad & &\bar{\varphi}(0)=A-iB,\\
\dot {\bar{\psi}} & = & -\bar{\varphi}(z\bar{\psi}+\bar{w} \varphi), \quad & &\bar{\psi}(0)=C-iD.
\end{array}
\end{equation*}
Here the constants $A,B,C,$ and $D$ have the following dynamical meaning: $\dot w(0)=C+iD$,
and $2B=i \ddot w(0)/\dot w(0)$. This complex  Hamiltonian system has the first integrals
\[
\begin{array}{rcl}
z\psi +w \varphi  & = &  C+iD, \\
z\bar{\varphi}+w \bar{\psi} & = & A-iB,\\
\end{array}
\]
and we have $|z|^2-|w|^2=1$ and $H=C^2+D^2=1$ as an additional normalization. Therefore, 
\[
\begin{array}{rcl}
\varphi  & = & z(A+iB)-\bar w ( C+iD), \\
\psi  & = & \bar z (C+iD) -w (A+iB).\\
\end{array}
\]
Let us introduce an auxiliary function $p=\bar w/{z}$. Then substituting $\varphi$ and $\psi$ in the Hamiltonian system we get
\[
p(s)=-(C-iD)\frac{1+e^{-2s\sqrt{1-B^2}}}{\sqrt{1-B^2}-iB+(\sqrt{1-B^2}-iB)e^{-2s\sqrt{1-B^2}}}.
\]
Taking into account that $\dot z \bar z= w \dot{\bar w}$, we get the solution for $B\neq 1$
\[
z(s)=\exp \int_0^s\frac{\bar{p}(t)\dot{p}(t)}{1-|p(t)|^2}dt,
\]
and
\[
w(s)=\bar{p}(s)\exp \int_0^s\frac{{p}(t)\dot{\bar p}(t)}{1-|p(t)|^2}dt.
\]
For $B=1$ the solution is
\[
z(s)=(1+is)e^{is},\quad w(s)=s(C+iD)e^{-is}.
\]
\medskip

\noindent
{\bf Parametric solution.}
Let us present the parametric solution in this case.
We use the parametrization in a neighbourhood of $(1,0,0,0)$ given by
\begin{equation}\label{e213} 
\begin{array}{ccr}
x_1 & = & \cos\chi_1\cosh\phi, \\ 
x_2 & = & \sin\chi_1\cosh\phi, \\ 
x_3 & = & \cos\chi_2\sinh\phi,\\
x_4 & = & \sin\chi_2\sinh\phi,
\end{array}
\end{equation} where $\phi\in(-\infty,\infty)$, $\chi_1,\chi_2\in(-\frac{\pi}{2},\frac{\pi}{2})$. We observe that the solution with constant velocity coordinates~\eqref{surface2} follows from this parameterization when we set $\phi=s$, $\chi_1=0$, and $\chi_2=\psi$. The vertical line (circle) $\Gamma$ is obtained by setting $\phi=0$, $\chi_1=s$, and $\chi_2=0$.

In this parametrization, the vector fields $T$, $X$, and $Y$ admit the form  
\begin{eqnarray*} 
T & = & \partial_{\chi_1}-\partial_{\chi_2} , \nonumber  \\ 
X & = & 2\cos(\chi_1-\chi_2)\partial_{\phi}-\partial_{\chi_1}\tanh\phi\sin(\chi_1-\chi_2)+\partial_{\chi_2}\cotanh\phi\sin(\chi_1-\chi_2), \nonumber \\
Y & = & 2\sin(\chi_1-\chi_2)\partial_{\phi}-\partial_{\chi_1}\tanh\phi\cos(\chi_1-\chi_2)+\partial_{\chi_2}\cotanh\phi\cos(\chi_1-\chi_2).
\end{eqnarray*} The vertical direction is given by the constant vector field $T$. 

The Hamiltonian $H$ associated with the operator $$\mathcal L=\frac{1}{2}(X^2+Y^2)=\frac{1}{2}(4\partial^2_{\phi}+\tanh^2\phi\partial^2_{\chi_1}+\cotanh^2\phi\partial^2_{\chi_2}-2\partial_{\chi_1}\partial_{\chi_2})$$ is given as $$H(\phi,\chi_1,\chi_2,\psi,\xi_1,\xi_2)=\frac{1}{2}(4\psi^2+\xi_1^2\tanh^2\phi+\xi_2^2\cotanh^2\phi-2\xi_1\xi_2),$$ where we set $\partial_{\phi}=\psi$, $\partial_{\chi_1}=\xi_1$, and $\partial_{\chi_2}=\xi_2$.

Description of geodesics is collected in the following theorem.
\medskip

\begin{theorem}
The geodesics starting from the point $\phi(0)=0$, $\chi_1(0)=0$, $\chi_2(0)=\chi_2^{(0)}$ with some $\dot\phi(0)$, a constant value of $\dot\chi_2$, and an arbitrary $\dot\chi_1(s)$ have the following equations.
\smallskip

\noindent If $4\dot\chi_2^2=\dot\phi^2(0)$ then
\begin{itemize}
\item[$\bullet$]{$\sinh\phi(s)=\pm |C|s$,}
 \item[$\bullet$]{$\chi_1(s)=\dot\chi_2 s-\frac{\dot\chi_2}{\dot\phi(0)}\arctan\dot\phi(0) s$,}
\item[$\bullet$] {$\chi_2(s)=-\dot\chi_2 s+\chi_2^{(0)}$.}
\end{itemize}

\smallskip

\noindent If $4\dot\chi_2^2>\dot\phi^2(0)$ then
\begin{itemize}
\item[$\bullet$]{$\sinh\phi(s)=\pm \sqrt{\frac{\dot\phi^2(0)}{4\dot\chi_2^2-\dot\phi^2(0)}}\sin\Big(s\sqrt{4\dot\chi_2^2-\dot\phi^2(0)}\Big)$,}
 \item[$\bullet$]{$\chi_1(s)=\dot\chi_2 s-\frac{\dot\chi_2}{2|\dot\chi_2|}\arctan\Big(\sqrt{\frac{4\dot\chi_2^2-\dot\phi^2(0)}{4\dot\chi_2^2}}\tan\big(s\sqrt{4\dot\chi_2^2-\dot\phi^2(0)}\big)\Big)$,}
\item[$\bullet$] {$\chi_2(s)=-\dot\chi_2 s+\chi_2^{(0)}$.}
\end{itemize}

\smallskip

\noindent If $4\dot\chi_2^2<\dot\phi^2(0)$ then
\begin{itemize}
\item[$\bullet$]{$\sinh\phi(s)=\pm \sqrt{\frac{\dot\phi^2(0)}{\dot\phi^2(0)-4\dot\chi_2^2}}\sinh\Big(s \sqrt{\dot\phi^2(0)-4\dot\chi_2^2}\Big)$,}
 \item[$\bullet$]{$\chi_1(s)=\dot\chi_2 s-\frac{\dot\chi_2}{2|\dot\chi_2|}\arctan\Big(\sqrt{\frac{\dot\phi^2(0)-4\dot\chi_2^2}{4\dot\chi_2^2}}\cotan\big(s\sqrt{4\dot\chi_2^2-\dot\phi^2(0)}\big)\Big)$,}
\item[$\bullet$] {$\chi_2(s)=-\dot\chi_2 s+\chi_2^{(0)}$.}
\end{itemize}
\end{theorem}

\section{Acknowledgment}

The paper was initiated when the authors visited the National Center for
Theoretical Sciences and National Tsing Hua University during
January 2007. They would like to express their profound gratitude
to Professors Jing Yu and Shu-Cheng Chang for their invitation and
for the warm hospitality extended to them during their stay in
Taiwan.

\end{document}